\documentclass[a4paper,11pt]{article}
\usepackage{amsfonts,amssymb,latexsym}
\usepackage[cp1251]{inputenc}
\usepackage[T2A]{fontenc}
\usepackage{graphicx}
\usepackage{amsmath}
\usepackage{longtable}
\usepackage{multicol}
\begin{document}
\Large

\newcommand{\Dp}{\Delta^+}
\newcommand{\De}{\Delta}
\newcommand{\gog}{{\frak g}}
\newcommand{\mog}{{\frak m}}
\newcommand{\ut}{{\frak u}{\frak t}}
\newcommand{\gl}{{\frak g}{\frak l}}
\newcommand{\pog}{{\frak p}}
\newcommand{\iog}{{\frak i}}
\newcommand{\DC}{{\cal D}}
\newcommand{\AC}{{\cal A}}
\newcommand{\GC}{{\cal G}}
\newcommand{\CC}{{\cal C}}
\newcommand{\SC}{{\cal S}}
\newcommand{\BC}{{\cal B}}
\newcommand{\PC}{{\cal P}}
\newcommand{\IC}{{\cal I}}
\newcommand{\JC}{{\cal J}}
\newcommand{\ogog}{\overline{\gog}}
\newcommand{\Nb}{{\Bbb  N}}
\newcommand{\ad}{{\mathrm{ad}}}
\newcommand{\eps}{\varepsilon}
\newcommand{\vphi}{\varphi}
\newcommand{\Ann}{\mathop{\mathrm{Ann}}\nolimits}
\newcommand{\Ker}{\mathop{\mathrm{Ker}}\nolimits}
\newcommand{\Maxspec}{\mathop{\mathrm{Maxspec}}\nolimits}
\newcommand{\lee}{\leqslant}
\newcommand{\gee}{\geqslant}
\newcommand{\Omegas}{\Omega_{\mathrm{sreg}}}

\author{M.V. Ignatev,\and A.N. Panov\footnote{The research was supported by RFFI grants
06-01-00037
 and 05-01-00313.}}
\date{}
\title{Coadjoint orbits of the group $UT(7,K)$}
 \maketitle

\section*{ \S0. Introduction}
Orbits of the coadjoint action play an important role in harmonic
analysis, theory of dynamic systems, noncommutative geometry. A.A.
Kirillov's orbit method allows to reduce the classification problem
of unitary representations of nilpotent Lie groups to the
classification of coadjoint orbits ~\cite{K-Orb,K-62}. This makes
possible to solve the problems of representation theory in geometric
terms of the orbit space. But it turns out that the classification
of coadjoint orbits is a difficult problem itself. In particular the
classification problem of coadjoint orbits of the unitriangular
group $\mathrm{UT}(n,K)$ for an arbitrary $n$ is far from its
solution~\cite{K-Orb, K-Var,K-M}. The classification of regular
orbits (i.e, orbits of maximum dimension) was achieved in the
pioneering paper on the orbit method~\cite{K-62}. The classification
of coadjoint orbits for $n\leq 6$ is known~\cite[\S 6.3.3]{K-Orb}.

Number of coadjoint orbits of given dimension over the finite field
${\Bbb F}_q$ is po\-ly\-no\-mi\-al in $q$. In~\cite{K-Var,K-M} a
correspondence between these polynomials and Euler-Bernoulli
$q$-polynomials is conjectured; this conjecture was checked by
computer for $n\leq 11$. Authors of the paper~\cite{G-Sh} made
computer program finding coadjoint orbits over the finite field (for
their investigations in completely integrable Toda lattices).

Every coadjoint orbit in the space $\ut^*(n,K)$ is contained in some
orbit of left-right action of ${\mathrm{UT}}(n,K)\times
{\mathrm{UT}(n,K)}$. In~\cite{C1, C2,C3} the classification of
orbits of left-right action is presented and correspondent
representations are studied.

In this paper we consider the case $n\leq 7$ and subregular orbits
of an arbitrary di\-men\-si\-on. The paper consists of three
sections. The first section contains the classification of coadjoint
orbits of the unitriangular group for $n\leq 7$ over an arbitrary
field $K$ of zero characteristic (Theorem~1.9). The description of
orbits is given in terms of admissible subsets; it's convenient to
represent them as diagrams. In the next Theorem 1.10 we solve the
classification problem of coadjoint orbits for $n\leq 7$ in terms of
canonical forms. For any canonical form the polarization is
constructed. This allows to classify unitary irreducible
representations of the real unitriangular groups and absolutely
maximal primitive ideals of the respective universal enveloping
algebras (Theorem 1.12 and Corollaries). Note that this diagram
method doesn't work for the unitriangular group of an arbitrary
size. Authors have constructed counterexamples to Theorems 1.3 and
1.10 for $n=9$. At the end of the first section we prove that the
number of orbits of Borel subgroup $\mathrm{T}(n,K)$ in $\ut^*(n,K)$
is infinite if $n>5$ (Theorem 1.15).

As usual, we identify the  symmetric algebra $S(\gog)$ with the
algebra $K[\gog^*]$ of regular functions on $\gog^*$. Coadjoint
orbits of the group $\mathrm{UT}(n,K)$ are closed, because orbits of
a regular action of a nilpotent group on an affine variety are
closed~(\cite[11.2.4]{Dix}). Let $I(\Omega)$ be the ideal of
functions that equal zero on the orbit $\Omega$ in $S(\gog)$. This
ideal is an absolutely maximal ideal (AMP-ideal) with respect to the
natural Poisson bracket in $S(\gog)$. The map $\Omega\mapsto
I(\Omega)$ establishes the bijection between the set of coadjoint
orbits and the set of AMP-ideals. In \S 2 we describe the set of
generators of an arbitrary AMP-ideal for $n\leq 7$ (Theorem 2.2).
This allows to represent an orbit as a set of solutions of
polynomial equations. It turns out that one can present the system
of generators as a system of some polynomials of the form $P-c$,
where $P$ is a coefficient of a minor of the characteristic matrix
$\Phi(\tau)$, $c\in K$. Authors believe that the last proposition
(see Corollary 2.4) is true for all coadjoint orbits of the
unitriangular group of an arbitrary dimension.

The description of regular orbits is given in terms of minors of the
matrix $\Phi$ (~\cite{K-62} and Theorem 3.1). In \S 3 of the paper
(Theorem 3.3) we present the classification of subregular orbits
(i.e., orbits of dimension $\dim \Omega_{\rm{reg.}}-2$). The
solution is given in terms of the coefficients of minors of the
characteristic matrix.

\section*{ \S1. Admissible diagrams and canonical forms of orbits}
Let $K$ be a field of zero characteristic. The unitriangular group
$\mathrm{UT}(n,K)$ consists of the matrices
$$
\left(\begin{array}{cccc}1&0&\ldots&0\\
a_{21}&1&\ldots&0\\
\vdots&\vdots&\ddots&\vdots\\
a_{n,1}&a_{n,2}&\ldots&1\end{array}\right),
$$
where $a_{ij}\in K$. The Lie algebra $\gog:=\ut(n,K)$ of the group
$G:=\mathrm{UT}(n,K)$ consists of lower-triangular matrices with
zeroes on the diagonal. Let $\Dp:=\Dp_n$ be the system of positive
roots of $\gl(n,K)$. All positive roots are
$\alpha_{ji}(h)=x_j-x_i$, where $j<i$ and
$h=\mathrm{diag}(x_1,\ldots,x_n)$.

For arbitrary pair $i>j$ denote by $y_{ij}$ the matrix from $\gog$
with 1 in the $(i,j)$-th entry and zeros elsewhere. The matrix
$y_{ij}$ is an eigenvector for $\ad_h$ with eigenvalue
$-\alpha_{ji}$. We say that the pair $(i,j)$ corresponds to the
positive root $\alpha_{ji}$. For a positive root $\xi =\alpha_{ji}$
we also denote by $y_\xi$ the matrix $y_{ij}$.

Matrices  $\{y_{ij},~i>j\}$ form a basis of the Lie algebra $\gog$.
By $\Phi$ we denote the following  formal matrix:
 $$
\Phi = \left(\begin{array}{cccc}0&0&\ldots&0\\
y_{21}&0&\ldots&0\\
\vdots&\vdots&\ddots&\vdots\\
y_{n1}&y_{n2}&\ldots&0\end{array}\right).$$

Let $\alpha > \beta$ be the lexicographic order on $\Dp$. I.e., if
pairs $(i, j)$ and $(s, t)$ correspond to the roots $\alpha$ and
$\beta$ resp., then $\alpha>\beta$ means that $j<t$ or $j=t$, $i>s$.
In other words, $\alpha>\beta$ if the entry $(i,j)$ of the $n\times
n$-matrix lies at the left of $(s,t)$ or in the same column, but
lower, than  $(s,t)$.

There exists the Killing form $(x,y) ={\rm{Tr}}\,xy$ on the Lie
algebra $\gl(n,K)$. The Killing form is non-degenerate over a field
of zero characteristic. This allows to identify $\gog^*$ with the
subspace of upper-triangular matrices with zeroes on the main
diagonal. We also identify symmetric algebra $S(\gog)$ with algebra
$K[\gog^*]$ of polynomials on $\gog^*$. There exists the natural
Poisson bracket on $S(\gog)$ defined by the formula $\{x,y\}=[x,y]$
for all $x,y\in\gog$.\\
{\bf Definition 1.1}.\\
1. An ideal $I$ in $S(\gog)$ is a Poisson ideal if
 $\{I,S(\gog)\}\subset I$.\\
2. We say that an ideal $I$ is an absolutely maximal Poisson ideal
(AMP-ideal) if the ideal $I\otimes_K K'$ is a maximal Poisson ideal
in $S(\gog\otimes_K K')$ for any extension $K'/K$ of the ground
field.

One can see that an ideal $I$ is Poisson if and only if it's
$G$-invariant. The map $I\mapsto \Ann I$ establishes the bijection
between the set of all AMP-ideals of $S(\gog)$ and the set of
coadjoint orbits in $\gog^*$.\\
{\bf Definition 1.2.}\\
1. We say that $A\subset \Dp$ is an additive subset if
$\alpha+\beta\in A$ whenever  $\alpha,\beta\in A$ and
$\alpha+\beta\in\Dp$.\\
2. Let $M\subset A$. If $\alpha\in A,\beta\in M$ and
$\alpha+\beta\in\Dp$ imply $\alpha+\beta\in M$, then we say that $M$
is a normal subset of $A$.

If $A$ is an additive subset of $\Dp$ and $M$ is a normal subset of
$A$, then $\gog_A:=\rm{span}\{y_\gamma:~\gamma\in A\}$ is a
subalgebra and $\mog:=\rm{span}\{y_\gamma:~\gamma\in M\}$ is an
ideal of $\gog_A$.

 For an arbitrary $\xi$ in an additive subset $A\subset\Dp$ consider
$$ C(\xi,A)= \{\gamma\in  A:~\xi>_A\gamma\}.$$

The subset $C(\xi,A)$ splits into the union
$C(\xi,A)=C_+(\xi,A)\sqcup C_-(\xi,A)$, where $C_+(\xi,A)$ (resp.
$C_-(\xi,A)$) consists of $\gamma\in C(\xi,A)$ such that $\gamma
> \xi-\gamma$ (resp. $\gamma < \xi-\gamma$).

Denote by $A(\xi)$ the subset $ A\setminus C(\xi,A)$. One can see
that $A(\xi)$ ia also an additive subset.
\\
{\bf Definition 1.3}. The subset $S=\{\xi_1>\ldots>\xi_k\}\subset
\Delta^+$ is called admissible if it's constructed by the following
rule. The first element $\xi_1$ is an arbitrary element of
$A_1:=\Dp$. Let $A_2:= A_1(\xi_1)$. The second element $\xi_2$ is an
arbitrary element of $A_2$. The next element $\xi_3$ is an arbitrary
element of $A_3:= A_2(\xi_2)$ etc.

Let  $A:=A(S):= A_{k+1}$ and $M:=M(S):=A\setminus S$. Note that $M$
is a normal subset of $A$.

One can split the subset $S$ into the union of two subsets
$S_\otimes\sqcup S_\square$. By definition, $\xi_i\in S_\otimes$, if
$ A_i\ne A_{i+1}$ (i.e., $C(\xi_i,A_i)\ne \emptyset$), and $\xi_i\in
S_\Box$, if $ A_i= A_{i+1}$ (i.e., $C(\xi_i,A_i)= \emptyset$).

To each admissible subset $S=\{\xi_1>\ldots>\xi_k\}$ we assign the
diagram $D(S)$ which is said to be admissible. The diagram $D(S)$ is
the $n\times n$-matrix. Its $(i, j)$-entries are empty for all
$i\leq j$; if $i>j$, then $(i, j)$-th entry is filled by one of the
symbols $\otimes$, $\square$, $\bullet$, $\pm$. In order to
construct the admissible diagram $D(S)$, consider the empty $n\times
n$-matrix $D$. Let $(i_1,j_1)$, $i_1>j_1$, corresponds to the first
root $\xi_1$ in $S$. If $i_1-j_1 >1$ (resp. $i_1-j_1 =1$), then put
the symbol $\otimes$ (resp. the symbol $\square$) into the
$(i_1,j_1)$-th entry of $D$. Put the symbol $\bullet$ into the
entries $(a,b)$, $b<j_1$, and $(a,j_1)$, $a>i_1$. Put also the
symbol $+$ into all entries $(a,j_1)$, $1<a<i_1$, lying in the
$j_1$-th column, and the symbol $-$ into all entries $(i_1,a)$,
$j_1<a<i_1$, lying in the $i_1$-th row.

Suppose that the entry $(i_2,j_2)$, $i_2>j_2$, corresponds to the
next positive root $\xi_2$ in $S$. Put the symbol $\bullet$ into all
empty entries $(a,b)$, $b<j_2$, and $(a,j_2)$, $a>i_2$. If there are
empty pairs of entries $(i_2,a)$ and $(a,j_2)$, where $j_2< a <
i_2$, in the "lower-triangle"\ part of our matrix (i.e., among the
entries with $i>j$), then put the symbol $\otimes$ into the entry
$(i_2,j_2)$ and the symbol $-$ (resp. $+$) into the entry $(i_2,a)$
(resp. $(a,j_2)$). If there are no empty pairs of entries in the
"lower-triangular"\ part of considering matrix, then fill the entry
$(i_2,j_2)$ by the symbol $\square$. Note that if one of the entries
$(i_2,a)$, $(a,j_2)$ is already  filled on the previous step, then
we don't fit the other. Then we apply this procedure to the roots
$\xi_3,\ldots,\xi_k$ in $S$. If, after all, there are an empty
entries in the "lower-triangular"\ part, then put the symbol
$\bullet$ into them. As result, we have the admissible diagram
$D(S)$.

Note that entries of $D(S)$, that corresponds to the roots in
$M(S)$, $S_\otimes$, $S_\square$, are filled by the symbols
$\bullet$, $\otimes$, $\square$ respectively. Recall that $A(S)=
M(S)\sqcup S_\otimes\sqcup S_\square$.
\\
{\bf Example}. To the admissible subset $S=\{\alpha_{13},
\alpha_{25}, \alpha_{35}, \alpha_{34}\}$ of the system of positive
roots for $n=5$ we assign the diagram
\begin{center} {\small
\begin{tabular}{|p{0.2cm}|p{0.2cm}|p{0.2cm}|p{0.2cm}|p{0.2cm}|}
\hline & & & &\rule{0pt}{0.5cm}\\
\hline + & & & &\rule{0pt}{0.5cm}\\
\hline $\otimes$ & $-$ & & &\rule{0pt}{0.5cm}\\
\hline $\bullet$ & + & $\square$ & &\rule{0pt}{0.5cm}\\
\hline $\bullet$ & $\otimes$ & $\square$ & $-$ &\rule{0pt}{0.5cm}\\
\hline\multicolumn{5}{c}{\rule{0pt}{0.5cm}}
\end{tabular}}
\end{center}
The construction of the diagram splits into four steps:
\begin{center} {\small
\begin{tabular}{|p{0.2cm}|p{0.2cm}|p{0.2cm}|p{0.2cm}|p{0.2cm}|}
\hline & & & &\rule{0pt}{0.5cm}\\
\hline + & & & &\rule{0pt}{0.5cm}\\
\hline $\otimes$ & $-$ & & &\rule{0pt}{0.5cm}\\
\hline $\bullet$ &  &  & &\rule{0pt}{0.5cm}\\
\hline $\bullet$ &&   &  &\rule{0pt}{0.5cm}\\
\hline\multicolumn{5}{c}{\rule{0pt}{0.5cm}}
\end{tabular}~$\Longrightarrow$
\begin{tabular}{|p{0.2cm}|p{0.2cm}|p{0.2cm}|p{0.2cm}|p{0.2cm}|}
\hline & & & &\rule{0pt}{0.5cm}\\
\hline + & & & &\rule{0pt}{0.5cm}\\
\hline $\otimes$ & $-$ & & &\rule{0pt}{0.5cm}\\
\hline $\bullet$ & + &  & &\rule{0pt}{0.5cm}\\
\hline $\bullet$ & $\otimes$ &  & $-$ &\rule{0pt}{0.5cm}\\
\hline\multicolumn{5}{c}{\rule{0pt}{0.5cm}}
\end{tabular}~$\Longrightarrow$
\begin{tabular}{|p{0.2cm}|p{0.2cm}|p{0.2cm}|p{0.2cm}|p{0.2cm}|}
\hline & & & &\rule{0pt}{0.5cm}\\
\hline + & & & &\rule{0pt}{0.5cm}\\
\hline $\otimes$ & $-$ & & &\rule{0pt}{0.5cm}\\
\hline $\bullet$ & + &  & &\rule{0pt}{0.5cm}\\
\hline $\bullet$ & $\otimes$ & $\square$ & $-$ &\rule{0pt}{0.5cm}\\
\hline\multicolumn{5}{c}{\rule{0pt}{0.5cm}}
\end{tabular}~$\Longrightarrow$
\begin{tabular}{|p{0.2cm}|p{0.2cm}|p{0.2cm}|p{0.2cm}|p{0.2cm}|}
\hline & & & &\rule{0pt}{0.5cm}\\
\hline + & & & &\rule{0pt}{0.5cm}\\
\hline $\otimes$ & $-$ & & &\rule{0pt}{0.5cm}\\
\hline $\bullet$ & + & $\square$ & &\rule{0pt}{0.5cm}\\
\hline $\bullet$ & $\otimes$ & $\square$ & $-$ &\rule{0pt}{0.5cm}\\
\hline\multicolumn{5}{c}{\rule{0pt}{0.5cm}}
\end{tabular}
}
\end{center}
Here $M(S)=\{\alpha_{15},\alpha_{14}\}$ and $A(S) =
\{\alpha_{15},\alpha_{14}\} \sqcup S$. Note that the subsets
$\{\alpha_{13}, \alpha_{25}, \alpha_{35}\}$, $\{\alpha_{13},
\alpha_{25},  \alpha_{34}\}$ and $\{\alpha_{13}, \alpha_{25}\}$ are
also admissible. Corresponding diagrams are
\begin{center} {\small
\begin{tabular}{|p{0.2cm}|p{0.2cm}|p{0.2cm}|p{0.2cm}|p{0.2cm}|}
\hline & & & &\rule{0pt}{0.5cm}\\
\hline + & & & &\rule{0pt}{0.5cm}\\
\hline $\otimes$ & $-$ & & &\rule{0pt}{0.5cm}\\
\hline $\bullet$ & + & $\bullet$ & &\rule{0pt}{0.5cm}\\
\hline $\bullet$ & $\otimes$ & $\square$ & $-$ &\rule{0pt}{0.5cm}\\
\hline\multicolumn{5}{c}{\rule{0pt}{0.5cm}}
\end{tabular}}
\quad\quad {\small
\begin{tabular}{|p{0.2cm}|p{0.2cm}|p{0.2cm}|p{0.2cm}|p{0.2cm}|}
\hline & & & &\rule{0pt}{0.5cm}\\
\hline + & & & &\rule{0pt}{0.5cm}\\
\hline $\otimes$ & $-$ & & &\rule{0pt}{0.5cm}\\
\hline $\bullet$ & + & $\square$ & &\rule{0pt}{0.5cm}\\
\hline $\bullet$ & $\otimes$ & $\bullet$ & $-$ &\rule{0pt}{0.5cm}\\
\hline\multicolumn{5}{c}{\rule{0pt}{0.5cm}}
\end{tabular}}
\quad\quad {\small
\begin{tabular}{|p{0.2cm}|p{0.2cm}|p{0.2cm}|p{0.2cm}|p{0.2cm}|}
\hline & & & &\rule{0pt}{0.5cm}\\
\hline + & & & &\rule{0pt}{0.5cm}\\
\hline $\otimes$ & $-$ & & &\rule{0pt}{0.5cm}\\
\hline $\bullet$ & + & $\bullet$ & &\rule{0pt}{0.5cm}\\
\hline $\bullet$ & $\otimes$ & $\bullet$ & $-$ &\rule{0pt}{0.5cm}\\
\hline\multicolumn{5}{c}{\rule{0pt}{0.5cm}}
\end{tabular}}
\end{center}

These diagrams have common $S_\otimes$-part and different
$S_\square$-parts.

One can see that if the subsets $S=S_\otimes \sqcup S_\Box$ and
$S'=S_\otimes \sqcup S'_{\square}$ are admissible then the subset
$S_\otimes \sqcup (S_\Box\cup S_\Box')$ is also admissible.
\\
 {\bf Definition 1.4}.
An admissible subset $S$ is called maximal if its subset $S_\Box$ is
maximal.

It's convenient to construct maximal admissible subsets for the
unitriangular Lie algebra of given dimension by the following
"sequence rule". The first maximal admissible subset is the
so-called regular subset $S_{\rm{reg}}$. Its first root $\xi_1$ is
the maximal root of $A_1:=\Dp$, the next root $\xi_2$ is the maximal
root in $A_2:= A_1(\xi_1)$ etc. (see the diagrams (n,0,1) for
$n=3,4,5,6,7$ above). Suppose that the maximal admissible subset
$S=\{\xi_1 > \ldots > \xi_p > \ldots > \xi_k\}$, where $\xi_p$ is
the minimal root in $S_\otimes$, is already constructed. Then the
next maximal admissible subset $S'=\{\xi'_1 > \ldots
>\xi'_{p}>\ldots
> \xi'_m\}$ is constructed as follows. The roots $\xi'_i$ coinside
with the roots $\xi_i$ for $1\leq i \leq p-1$, the root $\xi'_{p}$
is the maximal root of $\{\eta\in A_p:~ \xi_p>\eta\}$, and the root
$\xi'_i$ is the maximal root of $A_i$ for all $i>p$.

This procedure allows to list all maximal admissible diagrams. Here
we list all maximal admissible diagrams for $n=3,4,5$. Diagrams for
$n=6,7$ are listed at the end of the paper (we omit empty cells in
those diagrams). The number of the diagram consists of three numbers
$(n,k,m)$, where $n$ is the dimension of matrix, $k$ is the number
of symbols $\bullet$ in the first column and $m$ is the serial
number of the diagram in the set $(n,k)$. Every set $(n,k)$ ends by
the symbolic diagram $(n,k,\star)$, which corresponds to the set of
the diagrams: one can delete the first row and the first column and
replace the $(n-1)\times(n-1)$-block by an arbitrary $(n-1)\times
(n-1)$-diagram.

\bigskip{\small
\begin{center}
\begin{tabular}{cc}
\begin{tabular}{|p{0.2cm}|p{0.2cm}|p{0.2cm}|}
\hline & &\rule{0pt}{0.5cm}\\
\hline + & &\rule{0pt}{0.5cm}\\
\hline $\otimes$ & $-$ &\rule{0pt}{0.5cm}\\
\hline\multicolumn{3}{c}{\rule{0pt}{0.5cm}(3, 0, 1)}
\end{tabular}
&
\begin{tabular}{|p{0.2cm}|p{0.2cm}|p{0.2cm}|}
\hline & &\rule{0pt}{0.5cm}\\
\hline $\square$ & &\rule{0pt}{0.5cm}\\
\hline $\bullet$ & $\square$ &\rule{0pt}{0.5cm}\\
\hline\multicolumn{3}{c}{\rule{0pt}{0.5cm}(3, 1, 1)}
\end{tabular}
\end{tabular}

\par\bigskip
\begin{tabular}{ccc}
\begin{tabular}{|p{0.2cm}|p{0.2cm}|p{0.2cm}|p{0.2cm}|}
\hline & & &\rule{0pt}{0.5cm}\\
\hline + & & &\rule{0pt}{0.5cm}\\
\hline + & $\square$ & &\rule{0pt}{0.5cm}\\
\hline $\otimes$ & $-$ & $-$ &\rule{0pt}{0.5cm}\\
\hline\multicolumn{4}{c}{\rule{0pt}{0.5cm}(4, 0, 1)}
\end{tabular}

&
\begin{tabular}{|p{0.2cm}|p{0.2cm}|p{0.2cm}|p{0.2cm}|}
\hline & & &\rule{0pt}{0.5cm}\\
\hline + & & &\rule{0pt}{0.5cm}\\
\hline $\otimes$ & $-$ & &\rule{0pt}{0.5cm}\\
\hline $\bullet$ & $\square$ & $\square$ &\rule{0pt}{0.5cm}\\
\hline\multicolumn{4}{c}{\rule{0pt}{0.5cm}(4, 1, 1)}
\end{tabular}

&
\begin{tabular}{|p{0.2cm}|p{0.2cm}|p{0.2cm}|p{0.2cm}|}
\hline & & &\rule{0pt}{0.5cm}\\
\hline $\square$ & & &\rule{0pt}{0.5cm}\\
\hline $\bullet$ & * & &\rule{0pt}{0.5cm}\\
\hline $\bullet$ & * & * &\rule{0pt}{0.5cm}\\
\hline\multicolumn{4}{c}{\rule{0pt}{0.5cm}(4, 2, $\star$)}
\end{tabular}
\end{tabular}

\par\bigskip\begin{tabular}{cccc}
\begin{tabular}{|p{0.2cm}|p{0.2cm}|p{0.2cm}|p{0.2cm}|p{0.2cm}|}
\hline & & & &\rule{0pt}{0.5cm}\\
\hline + & & & &\rule{0pt}{0.5cm}\\
\hline + & + & & &\rule{0pt}{0.5cm}\\
\hline + & $\otimes$ & $-$ & &\rule{0pt}{0.5cm}\\
\hline $\otimes$ & $-$ & $-$ & $-$ &\rule{0pt}{0.5cm}\\
\hline\multicolumn{5}{c}{\rule{0pt}{0.5cm}(5, 0, 1)}
\end{tabular}

&
\begin{tabular}{|p{0.2cm}|p{0.2cm}|p{0.2cm}|p{0.2cm}|p{0.2cm}|}
\hline & & & &\rule{0pt}{0.5cm}\\
\hline + & & & &\rule{0pt}{0.5cm}\\
\hline + & $\square$ & & &\rule{0pt}{0.5cm}\\
\hline + & $\bullet$ & $\square$ & &\rule{0pt}{0.5cm}\\
\hline $\otimes$ & $-$ & $-$ & $-$ &\rule{0pt}{0.5cm}\\
\hline\multicolumn{5}{c}{\rule{0pt}{0.5cm}(5, 0, 2)}
\end{tabular}

&
\begin{tabular}{|p{0.2cm}|p{0.2cm}|p{0.2cm}|p{0.2cm}|p{0.2cm}|}
\hline & & & &\rule{0pt}{0.5cm}\\
\hline + & & & &\rule{0pt}{0.5cm}\\
\hline + & + & & &\rule{0pt}{0.5cm}\\
\hline $\otimes$ & $-$ & $-$ & &\rule{0pt}{0.5cm}\\
\hline $\bullet$ & $\otimes$ & $-$ & $\square$ &\rule{0pt}{0.5cm}\\
\hline\multicolumn{5}{c}{\rule{0pt}{0.5cm}(5, 1, 1)}
\end{tabular}

&
\begin{tabular}{|p{0.2cm}|p{0.2cm}|p{0.2cm}|p{0.2cm}|p{0.2cm}|}
\hline & & & &\rule{0pt}{0.5cm}\\
\hline + & & & &\rule{0pt}{0.5cm}\\
\hline + & $\square$ & & &\rule{0pt}{0.5cm}\\
\hline $\otimes$ & $-$ & $-$ & &\rule{0pt}{0.5cm}\\
\hline $\bullet$ & $\bullet$ & $\square$ & $\square$ &\rule{0pt}{0.5cm}\\
\hline\multicolumn{5}{c}{\rule{0pt}{0.5cm}(5, 1, 2)}
\end{tabular}

\par\bigskip \\
\begin{tabular}{|p{0.2cm}|p{0.2cm}|p{0.2cm}|p{0.2cm}|p{0.2cm}|}
\hline & & & &\rule{0pt}{0.5cm}\\
\hline + & & & &\rule{0pt}{0.5cm}\\
\hline $\otimes$ & $-$ & & &\rule{0pt}{0.5cm}\\
\hline $\bullet$ & + & $\square$ & &\rule{0pt}{0.5cm}\\
\hline $\bullet$ & $\otimes$ & $\square$ & $-$ &\rule{0pt}{0.5cm}\\
\hline\multicolumn{5}{c}{\rule{0pt}{0.5cm}(5, 2, 1)}
\end{tabular}

&
\begin{tabular}{|p{0.2cm}|p{0.2cm}|p{0.2cm}|p{0.2cm}|p{0.2cm}|}
\hline & & & &\rule{0pt}{0.5cm}\\
\hline + & & & &\rule{0pt}{0.5cm}\\
\hline $\otimes$ & $-$ & & &\rule{0pt}{0.5cm}\\
\hline $\bullet$ & $\square$ & + & &\rule{0pt}{0.5cm}\\
\hline $\bullet$ & $\bullet$ & $\otimes$ & $-$ &\rule{0pt}{0.5cm}\\
\hline\multicolumn{5}{c}{\rule{0pt}{0.5cm}(5, 2, 2)}
\end{tabular}

&
\begin{tabular}{|p{0.2cm}|p{0.2cm}|p{0.2cm}|p{0.2cm}|p{0.2cm}|}
\hline & & & &\rule{0pt}{0.5cm}\\
\hline + & & & &\rule{0pt}{0.5cm}\\
\hline $\otimes$ & $-$ & & &\rule{0pt}{0.5cm}\\
\hline $\bullet$ & $\square$ & $\square$ & &\rule{0pt}{0.5cm}\\
\hline $\bullet$ & $\bullet$ & $\bullet$ & $\square$ &\rule{0pt}{0.5cm}\\
\hline\multicolumn{5}{c}{\rule{0pt}{0.5cm}(5, 2, 3)}
\end{tabular}

&
\begin{tabular}{|p{0.2cm}|p{0.2cm}|p{0.2cm}|p{0.2cm}|p{0.2cm}|}
\hline & & & &\rule{0pt}{0.5cm}\\
\hline $\square$ & & & &\rule{0pt}{0.5cm}\\
\hline $\bullet$ & * & & &\rule{0pt}{0.5cm}\\
\hline $\bullet$ & * & * & &\rule{0pt}{0.5cm}\\
\hline $\bullet$ & * & * & * &\rule{0pt}{0.5cm}\\
\hline\multicolumn{5}{c}{\rule{0pt}{0.5cm}(5, 3, $\star$)}
\end{tabular}
\end{tabular}
\end{center}}

{\bf Definition 1.5}. Consider the set $\SC_n$, which consists of
all pairs $(S,c)$, where $S$ is a maximal admissible subset of
$\Dp:=\Dp_n$, and $c$ is a map $S\to K$ such that $c(\xi)\ne 0$ for
$\xi\in S_\otimes$. Any pair $(S,c)$ is called maximal admissible,
and the set $\SC_n$ is called the set of maximal admissible pairs.

Let the extension $c\colon A(S)\to K$ of $c$ be given by
$c(\gamma)=0$ for $\gamma\in M(S)$.\\
{\bf Definition 1.6}. Let $(S,c)$ be a maximal admissible pair. We
say that $f_{S,c}\in\gog^*$ is a canonical form if the following
condition hold:
$$f_{S,c}(y_\gamma) = \left\{ \begin{array}{l} c(\gamma),~~
{\mbox{если}}~~ \gamma\in A(S),\\
0, ~~ {\mbox{если}}~~  \gamma\in \Dp\setminus A(S).
\end{array} \right.$$
Note that the restriction of $f_{S,c}$ to $\gog_{A(S)}$ is a
character of $\gog_{A(S)}$.

The unitriangular Lie algebras form the sequence $$\gog=
\ut(n,K)\supset \ut(n-1,K)\supset\ldots\supset\ut(2,K),$$ where
$\ut(t,K)$ consists of matrices from $\ut(n,K)$ with zeroes in the
first $n-t$ columns. Systems of positive roots of the Lie algebras
$\{\gl(t,K)\}_{2\leq t\leq n}$ also form the sequence
$\Dp_n\supset\Dp_{n-1}\supset\ldots\supset\Dp_2$.

Denote $\De^{(n)}=\emptyset$ and $\De^{(t)}=\Dp_{n-t+1}\setminus
\Dp_{n-t}$. One can see that $\gamma\in \De^{(t)}$, if $y_\gamma$
lies in the $t$-th column of the matrix $\Phi$. The system of
positive roots splits into the union $$\Dp =
\De^{(1)}\sqcup\De^{(2)}\sqcup\ldots\sqcup\De^{(n-1)}.$$

We also put $ A^{(t)}:= \De^{(t)}\cap A$, $ M^{(t)}:= \De^{(t)}\cap
M$ $ S^{(t)}:= \De^{(t)}\cap S$, $ S_\otimes^{(t)}:= \De^{(t)}\cap
S_\otimes$, $ S_\square^{(t)}:= \De^{(t)}\cap S_\square$.

Let $m(1)=1$. We denote $$m(t)=\max\{ j ~\vert~~\xi_j\in
S^{(1)}\sqcup\ldots\sqcup S^{(t-1)}\} + 1,$$
$$B_t=\Dp_{n-t+1}\cap A_{m(t)}.$$
Note that $A^{(t)} = B_t\cap \De^{(t)}$. Subsets $\{B_t\}$ are
additive, they form the sequence $\Dp = B_1\supset
B_{2}\supset\ldots\supset B_{n-1}\supset B_n=\emptyset$.

To each admissible subset $S$ of $\Dp$ we assign the sequence of the
diagrams $\{D_t(S)\}_{1\le t\leq n}$ constructed as follows. The
diagram $D_1(S)$ (of dimension $n\times n$) coincides with $D(S)$.
One can construct the diagram $D_t(S)$ of dimension
${(n-t+1)\times(n-t+1)}$ by the following rule:\\
 1)~delete the first $t-1$ rows and $t-1$ columns,\\
 2)~omit the symbols (more precisely, the symbol $\{-\}$)
 corresponding to the weights, which don't contain in $B_t$.\\
{\bf Example}. $D(S)=(5, 2, 1).$\\

\begin{center}
{\small
\par\bigskip\begin{tabular}{cccc}
$D_1=\quad$\begin{tabular}{|p{0.2cm}|p{0.2cm}|p{0.2cm}|p{0.2cm}|p{0.2cm}|}
\hline & & & &\rule{0pt}{0.5cm}\\
\hline + & & & &\rule{0pt}{0.5cm}\\
\hline $\otimes$ & $-$ & & &\rule{0pt}{0.5cm}\\
\hline $\bullet$ & + & $\square$ & &\rule{0pt}{0.5cm}\\
\hline $\bullet$ & $\otimes$ & $\square$ & $-$ &\rule{0pt}{0.5cm}\\
\hline\multicolumn{5}{c}{\rule{0pt}{0.5cm}}
\end{tabular}

& $D_2=\quad$\begin{tabular}{|p{0.2cm}|p{0.2cm}|p{0.2cm}|p{0.2cm}|}
\hline & & &\rule{0pt}{0.5cm}\\
\hline & & &\rule{0pt}{0.5cm}\\
\hline + & $\square$ & &\rule{0pt}{0.5cm}\\
\hline $\otimes$ & $\square$ & $-$ &\rule{0pt}{0.5cm}\\
\hline\multicolumn{4}{c}{\rule{0pt}{0.5cm}}
\end{tabular}

& $D_3=\quad$\begin{tabular}{|p{0.2cm}|p{0.2cm}|p{0.2cm}|}
\hline & &\rule{0pt}{0.5cm}\\
\hline $\square$ & &\rule{0pt}{0.5cm}\\
\hline $\square$ & &\rule{0pt}{0.5cm}\\
\hline\multicolumn{3}{c}{\rule{0pt}{0.5cm}}
\end{tabular}

& $D_4=\quad$\begin{tabular}{|p{0.2cm}|p{0.2cm}|}
\hline &\rule{0pt}{0.5cm}\\
\hline &\rule{0pt}{0.5cm}\\
\hline\multicolumn{2}{c}{\rule{0pt}{0.5cm}}
\end{tabular}
\end{tabular}}
\end{center}

The subsets $B_1$, $B_2$, $B_3$, $B_4$ consist of the roots
corresponding to all non-empty entries of the diagrams $D_1$, $D_2$,
$D_3$, $D_4$ resp.

Denote by $\gog_t$ the subalgebra of $\ut(n-t+1, K)$ generated by
$\{y_\gamma:~\gamma\in B_t\}$. These subalgebras form the sequence
$$ \gog = \gog_1\supset
\ldots\supset\gog_t\supset\gog_{t+1}\supset\ldots\supset\gog_{n}=\{0\}.$$

Let $I_t$ be the ideal of $S(\gog_t)$ generated by $ y_\gamma -
c(\gamma)$, $\gamma\in A^{(t)}$. It's easy to check that $I_t$ is a
Poisson ideal of $S(\gog_t)$.

Denote by $\AC_m $ the Poisson algebra $K[p_1,\ldots,p_m;
q_1,\ldots, q_m]$, where $\{p_i,q_i\}=1$ and $ \{p_i,q_j\}=0$ for
$i\ne j$. Let $\AC_0=K$.

A Poisson algebra $\CC$ is called a tensor product of a Poisson
algebras $\BC_1$ and $\BC_2$ if $\CC=\BC_1\otimes\BC_2$ and $\{\BC_1,\BC_2\}=0$.\\
{\bf Lemma 1.7}. Let $n\leq 7$. For any $t$ there exists the Poisson
inclusion $$\theta_t: S(\gog_{t+1})\to S(\gog_t)/ I_t\eqno (1.1)$$
such that the Poisson algebra $S(\gog_t)/ I_t $ splits into the
tensor product of two Poisson subalgebras
$$S(\gog_t)/
I_t = \AC\otimes \theta_t S(\gog_{t+1}),\eqno (1.2)$$ where the
Poisson subalgebra $ \AC $ is isomorphic to $\AC_{k_t}$ for a
certain $k_t$.
\\
{\bf Proof}. Using the list of maximal admissible diagrams for
$n\leq 7$, we conclude that there are only three cases of filling of
a $t$-column: 1) the $t$-th column contains one symbol $\otimes$ and
the symbols $\square$ are situated below then the symbol $\otimes$
(or don't appear in this column); 2) the $t$-th column contains one
symbol $\otimes$ and one of the symbols $\square$ is situated above
then $\otimes$;
3) the $t$-th column doesn't contain the symbol $\otimes$. \\
{\bf Case 1}.~The $t$-th column contains one symbol $\otimes$ and
the symbols $\square$ are situated below then the symbol $\otimes$
(or don't appear in this column). Then $S_\otimes^{(t)}$ consists of
the one root; denote this root by $ \xi$.

Consider the subsets $C_+:=
C_+(\xi,B_t)=\{\gamma_1>\ldots>\gamma_k\}$ and $C_-:=
C_-(\xi,B_t)=\{\gamma'_1<\ldots<\gamma'_k\}$,  where
$\xi=\gamma_i+\gamma_i'$ for $1\leq i\leq m$. The subset $A^{(t)}$
is the union $$ A^{(t)} = \{\xi\}\sqcup M^{(t)}\sqcup
S_\square^{(t)}\sqcup C_+.$$ The subset $B_{t+1}$ has the form
$B_{t+1}=B_t\setminus (A^{(t)}\sqcup C_-)$.

Take the maximal element $\gamma_1\in C_+$ and the minimal element
$\gamma_1'\in C_-$ (recall that $\gamma_1+\gamma_1'=\xi$).

Let $B_t':=B_t\setminus \{\gamma_1'\}$. The subset $B_t'$ is
additive. Indeed, if $\gamma'_1=\eta+\eta'$, where $\eta,\eta'\in
A$, $\eta>\eta'$, then $\gamma_1+\eta\in C_+$. This contradicts the
choice of $\gamma_1$.

By $\gog_t'$ denote the subalgebra generated by $y_\eta$, $\eta\in
B_t'$. The elements $p=y_{\gamma_1'}$, $q=y_{\gamma_1}y_\xi^{-1}$
satisfy the condition $\{p,q\}=1$ and generate the subalgebra
isomorphic to $\AC_1$.

We claim that $\{q,a\}\equiv0\bmod I_t$ for an arbitrary $a\in
S(\gog_t')$. Indeed, for an arbitrary $\eta\in B_t'$ the root
$\gamma_1+\eta$ is contained in $ A^{(t)}$. Since
$\eta\ne\gamma_1'$, we conclude that $\gamma_1+\eta\neq\xi$. Since
$\gamma_1$ is a maximal element of $C_+$, we have
$\gamma_1+\eta\notin C_+$. A root $\xi_i$ from $S_\square$ can't be
equal to a sum of two roots from the correspondent $A_i$, so
$\gamma_1+\eta\notin S_\square$. Hence, $\gamma_1+\eta\in M^{(t)}$
and so $[y_{\gamma_1},y_\eta]=0\bmod I_t$ for an arbitrary $\eta\in
B_t'$. We conclude that $\{q,a\}=0\bmod I_t$ for all $a\in
S(\gog_t')$.

By definition, put
$$\widetilde{a}=\sum_{s=0}^\infty
\frac{(-1)^s}{s!}\ad_p^s(a)q^s=a-\ad_p(a)q+\frac{1}{2}\ad^2_p(a)q^2-\ldots,
\eqno(1.3)$$ where $a\in S(\gog_{t}')$ and $\ad_p(a)=\{p,a\}$ .

It's easy to check directly that $\{p,\widetilde{a}\}=
\{q,\widetilde{a}\}=0\bmod I_t$ and
$\{\widetilde{a},\widetilde{b}\}=\widetilde{\{a,b\}}\bmod I_t$
(\cite{Dix},4.7.5). The map $a\mapsto\widetilde{a}$ can be extended
to the morphism of Poisson algebras $\theta': S(\gog_{t}')\to
S(\gog_{t})/ I_t$. Denote by $I_t'$ the ideal of $S(\gog_t)$
generated by $I_t$ and $y_{\gamma_1}$. The ideal $I_t'$ coincides
with the kernel of $\theta'$. Extend $\theta'$ to the monomorphism
$$\theta': S(\gog_t')/I_t'\to
S(\gog_t)/I_t.$$ The Poisson algebra $S(\gog_{t})/I_t$ splits into
the tensor product $ \AC_1\otimes \rm{Im}\theta'$.

Then, for an arbitrary $1 < m \leq k$ we consider the additive
subset $B_t^{(m)}=B_t\setminus \{\gamma'_1,\ldots\gamma'_m\}$. By
$\gog_t^{(m)}$ denote the subalgebra generated by $y_\eta$, $\eta\in
B_t^{(m)}$, and by $I_t^{(m)}$ denote the ideal generated by
$I_t^{(m-1)}$ and $y_{\gamma_m}$. As in the case $m=1$, one can
construct the monomorphism $$ \theta^{(m)}:
S(\gog_t^{(m)})/I_t^{(m)}\to S(\gog_t)/I_t.$$ The Poisson algebra
$S(\gog_{t})/I_t$ splits into the tensor product $ \AC_m\otimes
\rm{Im}\theta^{(m)}$. For $m=k$ we have $$S(\gog_t^{(k)})/I_t^{(k)}=
S(\gog_{t+1})~~\mbox{and}~~
 S(\gog_{t})/I_t =  \AC_{k}\otimes \rm{Im}\theta^{(k)}.$$
{\bf Case 2}.~There are only two diagrams $D_t(S)$ such than one of
their $t$-th columns contains one symbol $\otimes$ and one of the
symbols $\square$ is situated above then $\otimes$:

\bigskip{\small
\begin{center}
\begin{longtable}{cc}

\begin{tabular}[b]{|p{0.2cm}|p{0.2cm}|p{0.2cm}|p{0.2cm}|}
\hline & & &\rule{0pt}{0.5cm}\\
\hline $\square$ & & &\rule{0pt}{0.5cm}\\
\hline + & & &\rule{0pt}{0.5cm}\\
\hline $\otimes$ & & $-$ &\rule{0pt}{0.5cm}\\
\hline\multicolumn{4}{c}{\rule{0pt}{0.34cm}}\\
\multicolumn{4}{c}{\rule{0pt}{0.5cm}$D_4(7, 3, 4)$}\\
\end{tabular}
\hspace{2cm}&\hspace{2cm}\begin{tabular}[b]{|p{0.2cm}|p{0.2cm}|p{0.2cm}|p{0.2cm}|}
\hline & & &\rule{0pt}{0.5cm}\\
\hline + & & &\rule{0pt}{0.5cm}\\
\hline $\square$ & $\square$ & &\rule{0pt}{0.5cm}\\
\hline $\otimes$ & $-$ & &\rule{0pt}{0.5cm}\\
\hline\multicolumn{4}{c}{\rule{0pt}{0.34cm}}\\
\multicolumn{4}{c}{\rule{0pt}{0.5cm}$D_4(7, 3, 8)$}
\end{tabular}

\end{longtable}
\end{center}}

For these diagrams we construct the required inclusion directly.
Here $\gog$ is the Lie subalgebra generated by $y_\eta$
corresponding
to the non-empty entries of diagrams. \\
1) $ D_4(7,3,4)$. The ideal $I$ is generated by the elements
$y_{41}-c_1, y_{21}-c_2$, where $c_1\ne 0$. The factor algebra
$S(\gog)/I$ coincides with the algebra  $\AC_1$ generated by
$p=y_{43}y_{41}^{-1}$ and
$q=y_{31}$.\\
2) $D_4(7,3,8)$. The ideal $I$ is generated by the elements $
y_{41}-c_1, y_{31}-c_2$ with $c_1\ne 0$. The factor algebra
$S(\gog)/I$ has the form $\AC_1\otimes K[\widetilde{y}_{32}]$, where
$\widetilde{y}_{32} = y_{32} - c_1^{-1}y_{42}y_{31}$ and $\AC_1$ is
generated by $p=y_{42}y_{41}^{-1}$ and
$q=y_{21}$.\\
{\bf Case 3}.~The $t$-th column doesn't contain the symbol
$\otimes$. Then $ A^{(t)} = M^{(t)} \sqcup S_\square^{(t)}$. The
natural inclusion $\gog_{t+1}\to \gog_t$ can be extended to the
isomorphism $S(\gog_{t+1}) \to S(\gog_t)/I_t$. $\Box$

Let $\gog_t$, $I_t$, $p_i,q_i$ be as in Lemma 1.7. Denote
$$X_t : = \Ann  I_t\subset \gog_t^*,$$
$$
X_{t,0} :=\{ f\in X_t:~p_i(f)=q_i(f)= 0 ~\mbox{для}~~ 1\leq i\leq
k_t\}. $$ Let $\theta_{t*}$ be the map $X_t\to\gog_{t+1}^*$ induced
by the inclusion $\theta_t$ (see (1.1)). Since $I_t$ is a Poisson
ideal, then $X_t$ is a Poisson submanifold and $\theta_{t*}$ is a
Poisson map.\\
{\bf Corollary 1.8}.  We claim that\\
1) the map $\omega\mapsto \theta_{t*}^{-1}(\omega)$ is a bijection
between the set of all symplectic leaves of $\gog_{t+1}^*$ and the
set of all symplectic leaves of $X_t$;\\
2) the restriction $$\rho_t:=\theta_{t*}\vert_{X_{t,0}}$$ is a
bijection between $ X_{t,0}$ and $\gog_{t+1}^*$; this restriction
coincides with the restriction of the natural projection
$\gog_t\to\gog_{t+1}$ to $ X_{t,0}$. \\
{\bf Proof}. An arbitrary Poisson ideal of the tensor product
$\AC_k\otimes \BC$ has the form $\AC_k\otimes \JC$, where $\JC$ is a
certain Poisson ideal of $\BC$. So, the first part of the corollary
follows from (1.2).

To prove the second one, we note that for an arbitrary $a\in
S(\gog_{t+1})$ and $f\in X_{t,0}$ we have $\widetilde{a}(f) = a(f)$.
Indeed, it follows from (1.3), that
$$\widetilde{a}(f) = \sum_{s=0}^\infty
\frac{(-1)^s}{s!}\ad_p^s(a)(f)q^s(f) = a(f).~\Box$$

To a maximal admissible pair $(S,c)$ we assign the AMP-ideal
$\IC_{S,c}$ of $S(\gog)$, which is constructed by induction on the
number of column $1\leq t\leq n-1$:
$$\IC^{(1)}\subset \IC^{(2)}\subset\ldots\subset
\IC^{(n-1)}=:\IC_{S,c}.$$

By definition, the ideal $\IC^{(1)}$ coincides with $I_1$. According
to Lemma 1.7, $S(\gog)/I_1\cong \AC_{k_1}\otimes
\theta_1(S(\gog_2))$. Suppose that the ideal $\IC^{(t)}$ such that
$$
S(\gog)/\IC^{(t)}= \AC_k\otimes \theta(S(\gog_{t+1})),\eqno(1.4)
$$
where $k = k_1+\ldots+k_{t}$ and $\theta = \theta_1\ldots\theta_t$,
is already constructed. Then the ideal $\IC^{(t+1)}$ contains
$\IC^{(t)}$ and is generated modulo $\IC^{(t)}$ by $\theta
(I_{t+1})$.

So,
$$S(\gog)/\IC^{(t+1)}= \AC_k\otimes
\theta(S(\gog_{t+1})/\IC^{(t+1)})=
\AC_k\otimes\theta\left(\AC_{k_{t+1}}\otimes\theta_{t+1}(S(\gog_{t+1}))\right)$$
$$
= \AC_{k+k_{t+1}}\otimes \theta\theta_{t+1}(S(\gog_{t+2})).$$ Denote
by $\IC_{S,c}$ the ideal $\IC^{(n-1)}$. Since $\gog_{n}=0$, then
$$S(\gog)/\IC_{S,c}= \AC_{k_1+\ldots+k_{n-1}}.\eqno(1.5)$$

Let $\Omega_{S,c} = \Ann \IC_{S,c}$. It follows from (1.5) that the
ideal $\IC_{S,c}$ is a AMP-ideal, so $\IC_{S,c}$ coincides with the
defining ideal $\IC(\Omega_{S,c})$ of the orbit (this ideal consists
of all elements of $S(\gog)$ such that their restrictions to the
orbit $\Omega_{S,c}$ are equal to zero).\\
{\bf Theorem 1.9}. Let $n\leq 7$ and $\gog = \ut(n,K)$.\\
1) The map $(S,c)\mapsto \IC_{S,c}$ is a bijection between the set
$\SC_n$ of all maximal admissible pairs and the set of all
AMP-ideals of $S(\gog)$.\\
2) The map $(S,c)\mapsto \Omega_{S,c}$ is a bijection between the
set $\SC_n$ of all maximal admissible pairs and the set of all
coadjoint orbits of the group $\mathrm{UT}(n,K)$.\\
3) $\dim(\Omega_f)$ is equal to the number of symbols $\pm$ in the
diagram $D(S)$.\\
{\bf Proof}. Part 3) follows from (1.5). Part 2) is a corollary of
1). Let's prove the first part of the theorem.

Let $\IC$ be an AMP-ideal of $S(\gog)$. If $\IC$ contains all
elements $ y_{n1}, \ldots, y_{21}$ of the first column of the matrix
$\Phi$, then denote by $\IC^{(1)}$ the ideal $\langle y_{n1},
\ldots, y_{21}\rangle$ and consider the second column.

If $\IC$ contains $y_{n1}, \ldots, y_{i+1,1}$ and doesn't contain
$y_{i1}$, then the ideal $\IC_0^{(1)}$ generated by $\langle y_{n1},
\ldots, y_{i+1,1}\rangle$ is a prime Poisson ideal of $ S(\gog)$.
The element $y_{i1}$ is a Casimir element modulo $\IC_0^{(1)}$.
Hence, the ideal $\IC$ contains a certain element of the form
$y_{i1} - c$, $c\in K$. Consequently, $\IC$ contains the ideal
$\IC^{(1)}$ generated by $\IC_0^{(1)}$ and $y_{i1} - c$. Denote by
$\xi_1$ the positive root that corresponds to the entry $(i,1)$.

According to Lemma 1.7, $S(\gog)/\IC^{(1)}= \AC_k\otimes
\theta_1S(\gog_2)$. Consider the elements
$$\theta(y_{n2}), \ldots, \theta_1(y_{32})$$ (they correspond to
the basic vectors from the second column of the matrix $\Phi$). In
the same way as for the first column, we can construct
$S=\{\xi_1,\ldots\}$ and $c: S\to K$, starting from the ideal
$\IC$.~$\Box$
\\
{\bf Theorem 1.10}. Let $n\leq 7$ and let $(S,c)$ be a maximal
admissible pair.  The canonical linear form $f_{S,c}$ is contained
in the orbit $\Omega_{S,c}$.\\
{\bf Proof}. Let $\{\gog_t\}$ be the sequence of subalgebras
corresponding to $S$ as above. Corollary 1.8(1) allows to construct
the set of coadjoint orbits (or, equivalently, symplectic leaves)
$\Omega_t\in\gog_t^*$ (of the coadjoint action of the groups
${\rm{Exp}}(\gog_t)$ resp.) such that for an arbitrary $t$ the orbit
$\Omega_t$ is contained in $X_t$ and $\Omega_t =
\theta_{t*}^{-1}(\Omega_{t+1})$.

It clearly follows from Corollary 1.8(2) that the orbit
$\Omega_{S,c}$ contains the element ${\rho^{-1}_1\circ\cdots\circ
\rho^{-1}_{n-1}(0)}$; this element coincides with $f_{S,c}$.~$\Box$ \\
{\bf Corollary 1.11}. The ideal $\IC_{S,c}$ is generated by the
elements $Q_\eta - Q_\eta(f_{S,c})$ with $\eta\in A(S)$ and
$Q_\eta(f_{S,c})= c(\eta)$. Here an arbitrary element $Q_\eta$ has
the form $y_\eta+ R_{>\eta}$, where $R_{>\eta}$ is a certain element
of the subalgebra generated by $\{y_\gamma,~ \gamma > \eta \}$, and
$R_{>\eta}(f_{S,c})=0$.\\
{\bf Proof}.  The ideal $\IC^{(t)}$ is generated modulo
$\IC^{(t-1)}$ by the elements $\theta_{t-1}(y_\eta)-c(\eta)$,
$\eta\in A^{(t)}$. It follows from (1.3) that there exists the
element $Q_\eta = y_\eta+ R_{>\eta}$ of $S(\gog)$ (where $R_{>\eta}$
is as above) that coincides with $\theta_{t-1}(y_\eta)$. Finally,
$$Q_\eta(f_{S,c})=(y_\eta+ R_{>\eta})(f_{S,c})=
y_\eta(f_{S,c})= c(\eta).\Box$$

Let $\GC$ be a Lie algebra and $f\in \GC^*$. A subalgebra $\pog$ is
called a polarization of $f$, if $\pog$ is a maximal isotropic
subspace with respect to the skew-symmetric bilinear form $f([x,y])$ on $\GC$. \\
{\bf Theorem 1.12}. Let $n\leq 7$ and $S$ be a maximal admissible
subset. Let $\pog_S$ be the linear subspace spanned by the basic
vectors $\{y_\eta\}$ such that the root $\eta$ corresponds to the
symbols $+$, $\otimes$, $\square$ or $\bullet$ in the diagram $D(S)$
(i.e., not to the symbol $-$). If $S$ doesn't correspond to the
diagram (7,3,8), then  $\pog_S$ is a polarization of all linear
forms $f_{S,c}$. To construct the polarization  $\pog_{(7,3,8)}$,
one should replace $y_{54}$ by $y_{75}$.\\
{\bf Proof} is by direct enumeration of all diagrams of dimension
$\leq 7$.~$\Box$

The classification of coadjoint orbits and construction of
polarizations allow to classify all unitary irreducible
representations~\cite{K-Orb} and all primitive ideals of the
universal enveloping algebra~\cite{Dix}.\\
{\bf Corollary 1.13.} An arbitrary unitary irreducible
representation of the group $\mathrm{UT}(n,{\Bbb R})$ with $n\leq 7$
is induced from a certain one-dimensional representation
$e^{if\ln}$, $f=f_{S,c}$, of the subgroup $P_S=\exp(\pog_S)$. The
maximal admissible pair $(S,c)$ is uniquely determined by the
representation.~$\Box$\\
{\bf Corollary 1.14.} An arbitrary absolutely primitive ideal of
$U(\ut(n,K))$ with $n\leq 7$ is induced from a certain ideal of
$U(\pog_S)\Ker f\vert_{\pog_S}$, where $ f=f_{S,c}$. The maximal
admissible pair $(S,c)$ is uniquely determined by the
ideal.~$\Box$\\
Let ${\mathrm{T}}(n,K)$ be the Borel subgroup of
${\mathrm{GL}}(n,K)$. The adjoint action of the group
${\mathrm{T}}(n,K)$ maps the Lie algebra of unitriangular matrices
$\gog:= \ut(n,K)$ to itself. This defines the action of the group
${\mathrm{T}}(n,K)$ on $\gog^*$. \\
{\bf Theorem 1.15}. The number of ${\mathrm{T}}(n,K)$-orbits in
$\gog^*$ is finite for $n\leq 5$ and infinite for $n>5$.\\
{\bf Proof}. The group ${\mathrm{T}}(n,K)$ maps a given orbit
$\Omega_{S,c}$ to a certain orbit $\Omega_{S,c'}$, which corresponds
to the same maximal admissible subset $S$ and the other admissible
map $c'$.

One can see that for $n\leq 5$ an arbitrary set $\Omega_S$ of
coadjoint orbits $\{\Omega_{S,c}\}$ (with fixed $S$ and arbitrary
admissible map $c\colon S\to K$) is a ${\mathrm{T}}(n,K)$-orbit. So,
for $n\leq 5$ the number of ${\mathrm{T}}(n,K)$-orbits in $\gog^*$
if finite.

Then, consider the diagram (6,3,4) for $n=6$:

\begin{center}{\small\begin{tabular}{|p{0.2cm}|p{0.2cm}|p{0.2cm}|p{0.2cm}|p{0.2cm}|p{0.2cm}
|}
\hline & & & & & \rule{0pt}{0.5cm}\\
\hline + & & & & &\rule{0pt}{0.5cm}\\
\hline $\otimes$ & $-$ & & & & \rule{0pt}{0.5cm}\\
\hline $\bullet$ & $+$ & $\square$  & & & \rule{0pt}{0.5cm}\\
\hline $\bullet$ & $\otimes$ & $\square$  & $-$ & & \rule{0pt}{0.5cm}\\
\hline $\bullet$ & $\bullet$ & $\bullet$ & $\square$ & $\square$ & \rule{0pt}{0.5cm}\\
\hline
\end{tabular}}
\end{center}

This diagram corresponds to the maximal admissible subset $S$ (it
consists of six elements). Codimension of $\Omega_{S,c}$ in
$\Omega_S$ is equal to six. But codimension of arbitrary coadjoint
${\mathrm{UT}}(n,K)$-orbit in the correspondent
${\mathrm{T}}(n,K)$-orbit is less or equal to five. Hence, the
number of ${\mathrm{T}}(n,K)$-orbits is infinite.

For arbitrary $n$ one can construct an example of infinite set of
${\mathrm{T}}(n,K)$-orbits, using the diagram of the following form:
its last six rows and six columns form the diagram (6,3,4) and all
other entries of this diagram are filled by the symbol~$\bullet$.

Note that an orbit from the set (6,3,4) has the following defining
equations (the the next section):
$$ y_{61}=y_{51}=y_{41}=y_{62}=y_{63}=0,~~y_{31}=c_1\ne 0,~~ y_{52}=c_2\ne 0,$$
$$\left\vert\begin{array}{cc}
y_{42}&y_{43}\\ y_{52}&y_{53}\end{array}\right\vert = c_3, \quad
y_{53}y_{31} + y_{52}y_{21} = c_4, ~~ y_{64} = c_5, ~~  y_{65} =
c_6.~\Box$$

\section*{ \S2. Absolutely maximal Poisson ideals}

In this section we'll present the system of defining equations of
arbitrary coadjoint orbit for $n\leq 7$. Consider the so-called
characteristic matrix for the formal matrix $\Phi$ (see~\S1):

$$
\Phi(\tau) = \tau\Phi + E = \left(\begin{array}{cccc}1&0&\ldots&0\\
\tau y_{21}&1&\ldots&0\\
\vdots&\vdots&\ddots&\vdots\\
\tau y_{n1}&\tau y_{n2}&\ldots&1\end{array}\right).
$$

Recall that we identify $S(\gog)$ with the algebra of polynomials on
$\gog^*$ (we also identify $\gog^*$ with the space of all
upper-triangular matrices with zeroes on the main diagonal).

Every minor of the (formal) characteristic matrix $\Phi(\tau)$ is a
polynomial on $\tau$ with the coefficients in the symmetric algebra
$S(\gog)= K[\gog^*]$. Values of these coefficients on an element
$F\in \gog^*$ coincides with the coefficients of the correspondent
minor of the usual characteristic matrix $\tau F+E$.

According to Theorem 1.9, to an arbitrary admissible pair $(S,c)$
one can assign the AMP-ideal $\IC_{S,c}$ and the orbit
$\Omega_{S,c}=\Ann \IC_{S,c}$. Dimension of the correspondent orbit
is equal to the number of the symbols $\pm$ in the diagram. The
ideal $\IC_{S,c}$ coincides with $\IC(\Omega_{S,c})$. To each root
$\eta\in A(S)$ one can assign the generator of the ideal
$\IC_{S,c}$.  In terms of the diagram $D(S)$, to each symbol
$\square$, $\otimes$, $\bullet$ one can assign the generator of
$\IC_{S,c}$. In this section we'll show that there are coefficients
of minors of the characteristic matrix $\Phi(\tau)$, which generate
the ideal $\IC_{S,c}$.

Let  $S_\otimes=\{\beta_1>\ldots>\beta_{k_\otimes}\}$, where
$k_\otimes=\vert S_\otimes \vert$. Denote $\eta=\eps_j-\eps_i\in
A(S)$. The $(i,j)$-th entry of the diagram $D(S)$ ($i>j$) is filled
by one of the symbols $\square$, $\otimes$,  $\bullet$. To the root
$\eta$ assign the permutation
$$w_\eta:=s_{\beta_1}\cdots s_{\beta_t} s_\eta,
~~\mbox{where}\quad
\beta_1>\cdots>\beta_t>\eta\gee\beta_{t+1}>\cdots>\beta_{k_\otimes}.$$

Consider the set of the columns $\Lambda:=\Lambda_j=\{1,\ldots,j\}$
and the set of the rows
$$w_\eta(\Lambda)=\mathrm{ord}\{w_\eta(1),\ldots,w_\eta(j)\}.$$
The minor $M^\Lambda_{w_\eta\Lambda}(\tau)$ of the matrix $\Phi(t)$
is polynomial on $\tau$:
$$ M^\Lambda_{w_\eta\Lambda}(\tau)=P_{q,\eta}\tau^q+\ldots P_{d,\eta}\tau^d,\quad
q<\ldots < d,$$ where degree $q$ (resp. $d$) of the least (resp. the
leading) term is equal to
$$ q=\vert \Lambda \setminus w_\eta\Lambda\vert = \vert w_\eta\Lambda\setminus
\Lambda\vert,$$
$$d=\sharp\{1\leq m\leq j~\vert~ i_m > m \}.$$
Let $\phi_1,\ldots, \phi_{n-1}$ be fundamental weights.\\
{\bf Lemma 2.1}. Let $n\leq 7$. We claim that for an arbitrary
$\eta=\eps_j-\eps_i\in A(S) $ the weight $(1-w_\eta)\phi_j$ can be
uniquely represented as a sum $\eta+\sum \beta$, where $\beta$ runs
over a certain subset
$H(S,\eta)\subset\{\beta_1,\ldots,\beta_t\}\subset S_\otimes$.
\\
{\bf Proof} is by direct enumeration of all diagrams for $n\leq 7$.
$\Box$
\\
Denote $h = h(S,\eta):=|H(S,\eta)|+1$. Note that $h$ is uniquely
determined by $\eta$ and $S$.\\
{\bf Theorem 2.2}. Let $n\leq 7$. An AMP-ideal $\IC_{S,c}$ of
$S(\gog)$ is generated by the elements $P_{h,\eta} - P_{h,\eta}^0\in
K$, где $P_{h,\eta}^0 = P_{h,\eta}(f_{S,c})$.
 \\
{\bf Proof.} One can prove this theorem for every maximal admissible
subset $S$. For example, consider the subset (7,2,7):\\

\begin{center}{\small\begin{tabular}{|p{0.2cm}|p{0.2cm}|p{0.2cm}|p{0.2cm}|p{0.2cm}|p{0.2cm}|p{0.2cm}|}
\hline & & & & & &\rule{0pt}{0.5cm}\\
\hline + & & & & & &\rule{0pt}{0.5cm}\\
\hline + & + & & & & &\rule{0pt}{0.5cm}\\
\hline + & $\otimes$ & $-$ & & & &\rule{0pt}{0.5cm}\\
\hline $\otimes$ & $-$ & $-$ & $-$ & & &\rule{0pt}{0.5cm}\\
\hline $\bullet$ & $\bullet$ & + & $\square$ & $\square$ & &\rule{0pt}{0.5cm}\\
\hline $\bullet$ & $\bullet$ & $\otimes$ & $\square$ & $\square$ & $-$ &\rule{0pt}{0.5cm}\\
\hline
\end{tabular}}
\end{center}
\par\bigskip
The additive subset $A(S)$ has the from $A(S) = S_\otimes\sqcup
S_\square\sqcup M(S)$, where
$$S_\otimes=\{\alpha_{15},\alpha_{24}, \alpha_{37}\},\quad
S_\square =\{\alpha_{46}, \alpha_{47}, \alpha_{56},
\alpha_{57}\},\quad
 M(S)= \{\alpha_{16}, \alpha_{17}, \alpha_{26}, \alpha_{27}\}.$$

Let $\IC$ be the ideal generated by the elements $P_{h,\eta} -
P_{h,\eta}^0$ with $\eta\in A(S)$. To prove that $\IC=\IC_{S,c}$,
it's enough to check that\\
1)~$\IC$ is a Poisson ideal,\\
2)~$\IC(f_{S,c}) = 0$,\\
3)~the ideal $\IC$ is generated by the elements $y_\eta+ T_{>\eta} -
c$, where $\eta\in A(S)$ and $T _{>\eta}$ is contained in the
subalgebra generated by $\{y_\gamma,~\gamma>\eta \}$.

Indeed, 1), 2) imply $\IC\subset\IC_{S,c}$.  By 3) the algebra
$S(\gog)/\IC$ is isomorphic to the algebra of polynomials
$K[y_\eta:~\eta\in \Dp\setminus A(S)]$. Therefore  the ideal  $\IC$
is prime and  $\mathrm{Dim}\IC = |\Dp\setminus A(S)| =\mathrm{dim}
\Omega_{S,c}$. We conclude $\IC=\IC_{S,c}$.

Let us show that  $\IC$ really satisfies  1), 2) and 3). Denote
$P_{ij}= P_{h,\eta}$, where  $\eta = \alpha_{ji}$. By definition,
the ideal  $\IC$ is generated by  $P_{ij}-P_{ij}^0$ where
$P_{ij}^0=P_{ij}(f_{S,c})$ and  $(i ,j)$ runs over pairs which are
filled by  symbols  $\otimes$, $\bullet$, $\square$ in the diagram.
We consider  the lexicographical order on the set on monomials of
$y_{ij}$ such that  $y_\alpha
> y_\eta$ if $\alpha < \eta$.
 Here we list the  polynomials  $P_{ij}$ starting from the leading    term:
$$ P_{71}=y_{71},\quad P_{61}=y_{61},\quad P_{51}=y_{51},\quad
 P_{51}=y_{51},\quad P_{72}=\left\vert\begin{array}{cc} y_{51}&
y_{52}\\y_{71}&y_{72}\end{array}\right\vert,$$
$$
P_{62}=\left\vert\begin{array}{cc} y_{51}&
y_{52}\\y_{61}&y_{62}\end{array}\right\vert,\quad
P_{42}=\left\vert\begin{array}{cc} y_{41}&
y_{42}\\y_{51}&y_{52}\end{array}\right\vert, \quad P_{73}=
\left\vert\begin{array}{ccc} y_{41}&
y_{42}&y_{43}\\y_{51}&y_{52}&y_{53}\\
y_{71}&y_{72}&y_{73}
\end{array}\right\vert
 $$

$$ P_{74}= y_{74}\left\vert\begin{array}{cc} y_{41}&
y_{42}\\y_{51}&y_{52}\end{array}\right\vert +
y_{73}\left\vert\begin{array}{cc} y_{31}&
y_{32}\\y_{51}&y_{52}\end{array}\right\vert = y_{74} P_{42} + \ldots
,$$
$$ P_{64}= \left\vert\begin{array}{cc} y_{41}&
y_{42}\\y_{51}&y_{52}\end{array}\right\vert \cdot
\left\vert\begin{array}{cc} y_{63}&
y_{64}\\y_{73}&y_{74}\end{array}\right\vert = -  y_{64}y_{73}P_{42}
+ \ldots $$
$$
P_{75} = y_{75}y_{51}+ y_{74}y_{41}+ y_{73}y_{31} = y_{75}y_{51} +
\ldots,$$
$$P_{65}= \left\vert\begin{array}{cc} y_{63}&
y_{65}\\y_{73}&y_{75}\end{array}\right\vert y_{51} +
\left\vert\begin{array}{cc} y_{63}&
y_{64}\\y_{73}&y_{74}\end{array}\right\vert y_{41} =
-y_{65}y_{73}y_{51}+\ldots.
$$

The ideal  $\IC$ contains the Poisson ideal  $\IC^{(1)}$ generated
by the elements  $y_{71}$, $y_{61}$, $y_{51}-y_{51}^0$. The elements
 $ P_{62}, P_{72}$ are Casimir elements modulo
$\IC^{(1)}$.
 Extend the ideal  $\IC^{(1)}$ to the ideal  $\JC:= \langle\IC^{(1)}, P_{72}, P_{62}\rangle$.

One can directly check that the above polynomials  $ P_{42}$,
$P_{73}$, $P_{74}$, $P_{64}$, $P_{75}$, $P_{65}$ are Casimir
elements in the factor algebra modulo  $\JC$. Therefore the ideal
$\IC$ is a Poisson ideal; this proves 1). The proof of  2) is
obvious.

Choose the new system of generators in  $\IC$ that consists of
$y_{71}$, $y_{61}$, $y_{51}-y_{51}^0$, $y_{73}-y_{73}^0$,
$Q_{ij}-Q_{ij}^0$ where  $(i,j)$ ranges over
$(7,4),~(6,4),~(7,5),~(6,5)$ and where  $Q_{ij}$ is defined similar
to $P_{ij}$ changing  $ y_{51}, y_{73}, P_{42}$ by $ y^0_{51},
y^0_{73}, P^0_{42}$. This proves  3).
~$\Box$ \\
 {\bf Corollary 2.3}. An orbit $\Omega_{S,c}$ has the defining
equations $P_{h,\eta}=\rm{const}(\eta)$, where $\eta$ ranges  over
positive roots that correspond to the symbols $\bullet$, $\otimes$
and $\square$ in the diagram $D(S)$.~$\Box$\\
{\bf Corollary 2.4}. A coadjoint orbit for $n\leq 7$ has defining
equations of the form $P-c$, where $P$ is a certain coefficient of a
minor of the characteristic matrix and $c\in K$.~$\Box$

\section*{ \S3. Subregular orbits}

In this section we'll describe all subregular orbits of the
unitriangular group (Theorem 3.3). Denote
$n_0=\left[\frac{n}{2}\right]$, $n_\otimes =
\left[\frac{n-1}{2}\right]$, $N=\frac{n(n-1)}{2}$. Note that
$n=n_0+n_\otimes+1$.

At first, we'll recall the description of regular orbits. One can
see that the minors $$P_j:=M^{1,\ldots,j}_{n-j+1,\ldots,n},\quad
1\leq j\leq n_0\eqno(3.1),$$ (top indices are numbers of columns,
bottom indices are numbers of rows) of the matrix $\Phi$
are Casimir elements of $S(\gog)$.\\
{\bf Theorem 3.1}~\cite{K-62}. The defining ideal
$\Omega_{\rm{reg}}$ of a regular orbit $\Omega_{\rm{reg}}$ is
generated by the elements $P_1-c_1, \ldots,P_{n_0}-c_{n_0}$, where
$c_j\in K$ (and $c_j\ne 0$ for $1\leq j\leq n_\otimes$).

We'll recall a proof of this theorem and show that regular orbits
correspond to regular maximal admissible subsets (see \S 1). Note
that $\vert S_{\rm{reg}}\vert=n_0$ and $\vert
S_{\rm{reg},\otimes}\vert = n_\otimes$.

Now we define the polynomials $Z_1,\ldots, Z_{n_\otimes}$, which are
needed for the sequel. Consider the set of the minors
$$P_j(\tau):= M^{1,\ldots,j}_{n-j+1,\ldots,n}(\tau),\quad
1\leq j\leq n-1,$$ of the characteristic matrix $\Phi(\tau)$. For
$1\leq j\leq n_0$ the polynomial $P_j(\tau)$ is equal to
$P_j\tau^j$, where $P_j$ is the correspondent minor of $\Phi$.

For $j > n_0$ the least term of the polynomial $P_j(\tau)$ equals
$P_j\tau^{n-j}$. Denote by $Z_{n-j}$ the coefficient of
$\tau^{n-j+1}$ in the polynomial $P_j(\tau)$, $j > n_0$. Note that
$1\leq n-j < n - n_0 = n_\otimes+1$, because $n_0 < j \leq n-1$. In
particular $Z_1 = y_{n,n-1}y_{n-1,1} + \ldots + y_{n,2}y_{2,1}$ is
the coefficient of $\tau^2$ in the polynomial $M^{1,\ldots,
n-1}_{2,\ldots, n}(\tau)$.

The following Proposition 3.2 is needed for the sequel. Consider the
de\-com\-po\-si\-ti\-on of the space $X:=\gog^*$ to the subsets
$X=X_0\sqcup X_1\sqcup\ldots\sqcup X_{n-1}$, where $$X_i=\{f\in X
\vert~ f(y_{n1})=\ldots =f(y_{n-i+1,1})=0,~ f(y_{n-i,1})\ne 0 \}. $$
Denote $d_i =\max\{\dim\Omega\vert~ \Omega\subset X_i\}.$ We say
that an orbit is $i$-regular, if $\Omega\subset X_i$ and
$\dim\Omega=d_i$.

A maximal admissible subset $S=\{\xi_1>\ldots>\xi_k\}$ is called
$i$-regular (denote it by $S(X_i)$) if $\xi_1=\alpha_{1,n-i}$ and an
arbitrary element $\xi_t$, $2\leq t\leq k$, is maximal element of
the correspondent $A_t$ (see Definition 1.3). Note that entries of
the diagram $D(X_i)$ corresponding to $S(X_i)$, which are filled by
the
symbols $\bullet$, are exactly $(n,1),\ldots, (n-i+1,1)$. \\
{\bf Proposition 3.2}. We have
$$d_i=\left\{\begin{array}{ll} N-(n_0+2i), & {\mbox{if}} ~ i\leq
n_\otimes,\\ N-(n_0+2n_\otimes), & {\mbox{if}} ~ i>n_\otimes.
\end{array}\right.$$
{\bf Proof}. To each $S(X_i)$ we'll assign the dense in $X_i$ set of
orbits $\{\Omega_{S(X_i),c}\}$ of equal dimension. On the other
hand, the set of orbits of maximum dimension is open~\cite[\S
2.6]{Kraft}. Hence, $d_i = \dim \Omega_{S(X_i),c}$.

This set of orbits $\{\Omega_{S(X_i),c}\}$ will be constructed
differently in different cases 1) and 2).\\
{\bf Case 1.}~$i\leq n_\otimes$. See diagrams $D(X_2)$ for $n=7$ and
$n=8$.

{\small
\begin{center}
\begin{tabular}{ccc}
\begin{tabular}[b]{|p{0.2cm}|p{0.2cm}|p{0.2cm}|p{0.2cm}|p{0.2cm}|p{0.2cm}|p{0.2cm}|}
\hline\multicolumn{7}{|c|}{\rule{0pt}{0.5cm}}\\
\cline{1-1}  + &\multicolumn{6}{|c|}{\rule{0pt}{0.5cm}}\\
\cline{1-2}  + & + &\multicolumn{5}{|c|}{\rule{0pt}{0.5cm}}\\
\cline{1-3}  + & + & + &\multicolumn{4}{|c|}{\rule{0pt}{0.5cm}}\\
\cline{1-4}  $\otimes$ & $-$ & $-$ & $-$ &\multicolumn{3}{|c|}{\rule{0pt}{0.5cm}}\\
\cline{1-5}  $\bullet$ & + & $\otimes$ & $-$ & $\square$ &\multicolumn{2}{|c|}{\rule{0pt}{0.5cm}}\\
\cline{1-6}  $\bullet$ & $\otimes$ & $-$ & $-$ & $\square$ & $-$ &\rule{15pt}{0pt}\rule{0pt}{0.5cm}\\
\hline
\end{tabular}

&\quad \quad\quad\quad\quad&

\begin{tabular}[b]{|p{0.2cm}|p{0.2cm}|p{0.2cm}|p{0.2cm}|p{0.2cm}|p{0.2cm}|p{0.2cm}|p{0.2cm}|}
\hline\multicolumn{8}{|c|}{\rule{0pt}{0.5cm}}\\
\cline{1-1}  + &\multicolumn{7}{|c|}{\rule{0pt}{0.5cm}}\\
\cline{1-2}  + & + &\multicolumn{6}{|c|}{\rule{0pt}{0.5cm}}\\
\cline{1-3}  + & + & + &\multicolumn{5}{|c|}{\rule{0pt}{0.5cm}}\\
\cline{1-4}  + & + & + & $\square$ &\multicolumn{4}{|c|}{\rule{0pt}{0.5cm}}\\
\cline{1-5}  $\otimes$ & $-$ & $-$ & $-$ & $-$ &\multicolumn{3}{|c|}{\rule{0pt}{0.5cm}}\\
\cline{1-6}  $\bullet$ & + & $\otimes$ & $-$ & $-$ & $\square$ &\multicolumn{2}{|c|}{\rule{0pt}{0.5cm}}\\
\cline{1-7}  $\bullet$ & $\otimes$ & $-$ & $-$ & $-$ & $\square$ & $-$ &\rule{15pt}{0pt}\rule{0pt}{0.5cm}\\
\hline
\end{tabular}
\end{tabular}
\end{center}}

Columns of $D(X_i)$ satisfy conditions from the cases 1) or 3) from
the proof of Lemma 1.7. This allows to construct the ideal
$\IC_{S(X_i),c}$ by methods of Theorem 1.9. The dimension of the
orbit $\Omega_{S(X_i),c}$ is equal to the number of the symbols
$\pm$ in the diagrams. All orbits from this set have equal
dimension, which equals $ N - \vert S(X_i)\vert = N-(n_0+2i).$ This
concludes the proof in the fist case.

Note that one can find the ideal $\IC_{S(X_i),c}$ using methods of
\S 2. Indeed, let $\{\xi_1>\ldots>\xi_{n_0}\}$ be the first $n_0$
roots of $S(X_i)$ (i.e. the weights corresponding to the symbols
$\otimes$ in the odd case and the weights corresponding to the
symbols $\otimes$ with the weight $\alpha_{n_0,n_0+1}$ in the even
case).

To each weight $ 1\leq j\leq n_0$ assign the row number $m(j)$ such
that $\xi_i=\alpha_{i,m(i)}$. Consider the system of the minors
$Q_1,\ldots, Q_{n_0}$, where

$$Q_j = M^{1,\ldots, j}_{m(1),\ldots, m(j)}.$$

We also consider the ideal $\JC$ generated by $y_{n1},\ldots,
y_{n-i+1,1}$. It's easy to check that the polynomials $Q_1,\ldots,
Q_{n_0}, Z_1, \ldots, Z_{i}$ are Casimir elements modulo $\JC$. The
ideal $\IC_{S(X_i),c}$ is generated by $\JC$ and the elements of
the form $P-c$, where $P=Q_j$, $1\leq j\leq n_0$, or  $P=Z_j$, $1\leq j\leq i$.\\
{\bf Case 2.}~$i > n_\otimes$. See diagrams $D(X_4)$ for $n=7$ and
$n=8$:

{\small
\begin{center}
\begin{tabular}{ccc}
\begin{tabular}[b]{|p{0.2cm}|p{0.2cm}|p{0.2cm}|p{0.2cm}|p{0.2cm}|p{0.2cm}|p{0.2cm}|}
\hline\multicolumn{7}{|c|}{\rule{0pt}{0.5cm}}\\
\cline{1-1}  + &\multicolumn{6}{|c|}{\rule{0pt}{0.5cm}}\\
\cline{1-2}  $\otimes$ & $-$ &\multicolumn{5}{|c|}{\rule{0pt}{0.5cm}}\\
\cline{1-3}  $\bullet$ & + & + &\multicolumn{4}{|c|}{\rule{0pt}{0.5cm}}\\
\cline{1-4}  $\bullet$ & + & + & $\square$ &\multicolumn{3}{|c|}{\rule{0pt}{0.5cm}}\\
\cline{1-5}  $\bullet$ & + & $\otimes$ & $-$ & $-$ &\multicolumn{2}{|c|}{\rule{0pt}{0.5cm}}\\
\cline{1-6}  $\bullet$ & $\otimes$ & $\square$ & $-$ & $-$ & $-$ &\rule{15pt}{0pt}\rule{0pt}{0.5cm}\\
\hline
\end{tabular}
&\quad \quad\quad\quad\quad&
\begin{tabular}[b]{|p{0.2cm}|p{0.2cm}|p{0.2cm}|p{0.2cm}|p{0.2cm}|p{0.2cm}|p{0.2cm}|p{0.2cm}|}
\hline\multicolumn{8}{|c|}{\rule{0pt}{0.5cm}}\\
\cline{1-1}  + &\multicolumn{7}{|c|}{\rule{0pt}{0.5cm}}\\
\cline{1-2}  + & + &\multicolumn{6}{|c|}{\rule{0pt}{0.5cm}}\\
\cline{1-3}  $\otimes$ & $-$ & $-$ &\multicolumn{5}{|c|}{\rule{0pt}{0.5cm}}\\
\cline{1-4}  $\bullet$ & + & + & + &\multicolumn{4}{|c|}{\rule{0pt}{0.5cm}}\\
\cline{1-5}  $\bullet$ & + & + & $\otimes$ & $-$ &\multicolumn{3}{|c|}{\rule{0pt}{0.5cm}}\\
\cline{1-6}  $\bullet$ & + & $\otimes$ & $\square$ & $-$ & $-$ &\multicolumn{2}{|c|}{\rule{0pt}{0.5cm}}\\
\cline{1-7}  $\bullet$ & $\otimes$ & $-$ & $\square$ & $-$ & $-$ & $-$ &\rule{15pt}{0pt}\rule{0pt}{0.5cm}\\
\hline
\end{tabular}
\end{tabular}
\end{center}}

Columns of $D(X_i)$ also satisfy conditions from the cases 1) or 3)
from the proof of Lemma 1.7, so we can construct the ideal
$\IC_{S(X_i),c}$ by methods of Theorem 1.9. The dimension of
$\Omega_{S(X_i),c}$ equals the number of the symbols $\pm$ in the
diagrams. All orbits from this set have equal dimension, which is
equal to $ N - \vert S(X_i)\vert = N-(n_0+2n_\otimes).$ This
concludes the proof in the second case.

Note that one can find the generators of $\IC_{S(X_i),c}$ as in the
first case: the ideal $\IC_{S(X_i),c}$ is generated by $\JC$ and the
elements of the form $P-c$, where $P=Q_j$, $1\leq j\leq n_\otimes +
1$, or $P=Z_j$, $1\leq j\leq n-i-2$.~$\Box$\\
 {\bf Remark}. We've just shown that an orbit $\Omega_{S(X_i),c}$ is
$i$-regular. But if $i>1$, then an $i$-regular orbit
doesn't have the form $\Omega_{S(X_i),c}$ in general.\\
{\bf Proof of Theorem 3.1}. We can apply methods of Theorem 1.9 to
regular diagrams (i.e., diagrams $ D(S_{\rm{reg}})$); i.e., we can
construct the set of orbits $\{\Omega_{S_{\rm{reg}},c}\}$ of
dimension $N-n_0$, starting from $S_{\rm{reg}}$. Since this set is
dense in $\gog^*$, we conclude that the maximum dimension of an
orbit is equal to $N-n_0$. The ideal $\IC(S_{\rm{reg}})$ is
generated by $P_j-c_j$, $1\leq j\leq n_0$.

So, it's enough to check that an arbitrary regular orbit coincides
with one of the orbits constructed by a regular diagram. Let
$\Omega$ be a regular orbit. It follows from Proposition 3.2 that
$\Omega\subset X_0$. Since $y_{n1}$ are Casimir elements of
$S(\gog)$, that $\IC(\Omega)$ contains a certain element $y_{n1}-c$
with $c\ne 0$. Let $\IC^{(1)} = \langle y_{n1} -c\rangle$. Using
Lemma 1.7, we obtain
$$S(\gog)/\IC^{(1)} = \AC_{n-2} \otimes S(\ut(n-2,K)).$$
There exists the symplectic leaf (orbit) $\omega$ in $\ut^*(n-2,K)$,
corresponding to the symplectic leaf (orbit) $\Omega$ in $\gog^*$.
Since the orbit $\Omega$ is regular, then the orbit $\omega$ is
regular in $\ut^*(n-2,K)$. To conclude the proof, it remains to
apply induction on $n$.~$\Box$

Let's now describe subregular orbits. For $n=3$ and $n=4$ they are
described by the diagrams:

\bigskip{\small
\begin{center}
\begin{tabular}{ccc}

\begin{tabular}{|p{0.2cm}|p{0.2cm}|p{0.2cm}|}
\hline & &\rule{0pt}{0.5cm}\\
\hline $\square$ & &\rule{0pt}{0.5cm}\\
\hline $\bullet$ & $\square$ &\rule{0pt}{0.5cm}\\
\hline \multicolumn{3}{c}{\rule{0pt}{0.5cm}(3, 1, 1)}

\end{tabular}
\quad \quad & \quad\quad
\begin{tabular}{|p{0.2cm}|p{0.2cm}|p{0.2cm}|p{0.2cm}|}
\hline & & &\rule{0pt}{0.5cm}\\
\hline + & & &\rule{0pt}{0.5cm}\\
\hline $\otimes$ & $-$ & &\rule{0pt}{0.5cm}\\
\hline $\bullet$ & $\square$ & $\square$ &\rule{0pt}{0.5cm}\\
\hline \multicolumn{3}{c}{\rule{0pt}{0.5cm}(4, 1, 1)}

\end{tabular}
\quad \quad & \quad\quad
\begin{tabular}{|p{0.2cm}|p{0.2cm}|p{0.2cm}|p{0.2cm}|}
\hline & & &\rule{0pt}{0.5cm}\\
\hline $\square$ & & &\rule{0pt}{0.5cm}\\
\hline $\bullet$  & $+$ & &\rule{0pt}{0.5cm}\\
\hline $\bullet$ & $\otimes$ & $-$ &\rule{0pt}{0.5cm}\\
\hline \multicolumn{3}{c}{\rule{0pt}{0.5cm}(4, 2, 1)}

\end{tabular}
\end{tabular}
\end{center}}
One can check that for $n=5,6,7$ diagrams corresponding to
subregular orbits are: (5,0,2), (5,1,1) for $n=5$; (6,0,2), (6,0,3),
(6,1,1) for $n=6$; (7,0,2), (7,0,3), (7,1,1) for $n=7$.

For an arbitrary $1\leq j\leq n_\otimes$ consider the minors
$$P_j'=M^{1,\ldots,j}_{n-j, n-j+2,\ldots,n}\quad
 P_j''=M^{1,\ldots,j-1,j+1}_{n-j+1,\ldots,n}$$
of the matrix $\Phi$, which border the minor $P_{j-1}$ from (3.1).
For even $n$ consider also the minor $P_{n_0}' =
M^{1,2,\ldots,n_0-1}_{n_0, n_0+3, n_0+4,\ldots, n}$.

Recall that the restriction of a minor $P_j$ (from (3.1)) to an
arbitrary orbit is constant, because these minors are Casimir
elements.\\
{\bf Theorem 3.3}. Let $\Omega_{\rm{sreg}}$ be a subregular orbit.\\
1)~If $P_{n_\otimes}(\Omega_{\rm{sreg}})\ne 0$, then there are the
number $1\leq j_0 < n_\otimes$ and the numbers
$$\{c_1,\ldots, c_{j_0-1},c',c'',c_{j_0+2},\ldots, c_{n_0-1}\}\subset K^*,\quad
 \{c_{n_0}, c\}\subset K$$
(with $c_{n_0}\neq 0$ for odd $n$) such that the ideal
$I(\Omega_{\rm{sreg}})$ is generated by the elements $$P_i
-c_i,~\mbox{where}~i=1,\ldots, j_0-1, j_0+2,\ldots, n_0,$$
$$P_{j_0}'-c',\quad P_{j_0}''-c'',\quad P_{j_0},\quad Z_{j_0}-c.$$
2)~If $P_{n_{\otimes}}(\Omega_{\rm{sreg}}) = 0$ and $n$ is odd, then
there are the numbers
$$\{c_1,\ldots,c_{n_\otimes-1}\}\subset K^*,\quad \{c',c''\}\subset K$$
such that the ideal $I(\Omega_{\rm{sreg}})$ is generated by the
elements $$P_i-c_i, ~\mbox{where}~i=1,\ldots,n_{\otimes}-1,$$
$$P_{n_{\otimes}}'-c',\quad P_{n_{\otimes}}''-c'',\quad P_{n_{\otimes}}.$$
3)~If $P_{n_{\otimes}}(\Omega_{\rm{sreg}}) = 0$ and $n$ is even,
then there are the numbers
$$\{c_1,\ldots,c_{n_{\otimes}-1},c'\}\subset K^*,\quad
\{c'',c\}\subset K$$ such that the ideal $I(\Omega_{\rm{sreg}})$ is
generated by the elements $$P_i-c_i,
~\mbox{where}~i=1,\ldots,n_{\otimes}-1,$$ and the elements of the
form $a)$ or $b)$, where

$$\begin{array}{l} a) ~ P_{n_{\otimes}}'-c',~ P_{n_{\otimes}}''-c'',~ P_{n_{\otimes}},
~ Z_{n_{\otimes}}-c,~~ \mbox{и} \\
b) ~ P_{n_{\otimes}},~ P_{n_{\otimes}}',~ P_{n_0}'-c'',~
P_{n_{\otimes}}''-c'.\end{array}$$ {\bf Proof}. The proof for
$n\leq4$ is by direct calculations. Let $\Omega_{\mathrm{sreg}.}$ be
a subregular orbit in $\ut(n,K)$ for $n>4$. It follows from
Proposition 3.2 that there are two cases for $n>4$:
$\Omega_{\mathrm{sreg}}\subset X_0$ or
$\Omega_{\mathrm{sreg}}\subset X_1$.  \\
{\bf Cases 1.} $\Omega_{\mathrm{sreg}}\subset X_1$.

The ideal $\IC(\Omegas)$ contains the ideal $\JC$ generated by
$y_{n1}$. The elements $P_i$, $3\leq i\leq n_0$, and $y_{n-1,1}$,
$y_{n2}$, $Z_1$ are Casimir elements modulo $\JC$. There are $\{c,
c_1',c_1'', c_3\ldots, c_{n_0}\}\subset K$ such that the
restrictions to $\Omegas$ of all functions from the ideal $\IC$
generated by $\JC$, $P_i-c_i$, $3\leq i\leq n_0$ and $y_{n-1,1}-c'$,
$y_{n2} - c''$, $Z_1 - c$ with $c'\ne 0$, are equal to zero. Hence,
$\IC\subset \IC(\Omegas)$.

Consider the Poisson ideal $\IC_0 = \langle y_{n,1}, y_{n-1,1} - c',
Z_1 - c\rangle$. Easy to see that $\IC_0\subset \IC\subset
\IC(\Omegas)$. According to Lemma 1.7,
$$ S(\gog)/\IC_0 = \AC_{n-3}\otimes S(\gog_2),$$
where the Lie algebra $\gog_2$ is isomorphic to $\ut(n-2,K)$. Both
of ideals $\IC$, $\IC(\Omegas)$ are generated modulo $\IC_0$ by a
regular AMP-ideal of $S(\gog_2)$. Thus, $\IC =  \IC(\Omegas)$.
\\
{\bf Case 2.} $\Omega_{\mathrm{sreg}}\subset X_0$.

Suppose that $P_1(\Omegas)\ne 0,\ldots , P_{j_0-1}(\Omegas)\ne 0$
and $ P_{j_0}(\Omegas)=0$. The proof is similar to the first case
(one should replace $y_{n1}$, $y_{n-1,1}$, $y_{n2}$ by $P_{j_0}$,
$P'_{j_0}$, $P''_{j_0}$).~$\Box$

\newpage
\begin{center}
\bf{List of maximal admissible diagrams for $n=6$}
\par\bigskip
\includegraphics[width=16cm]{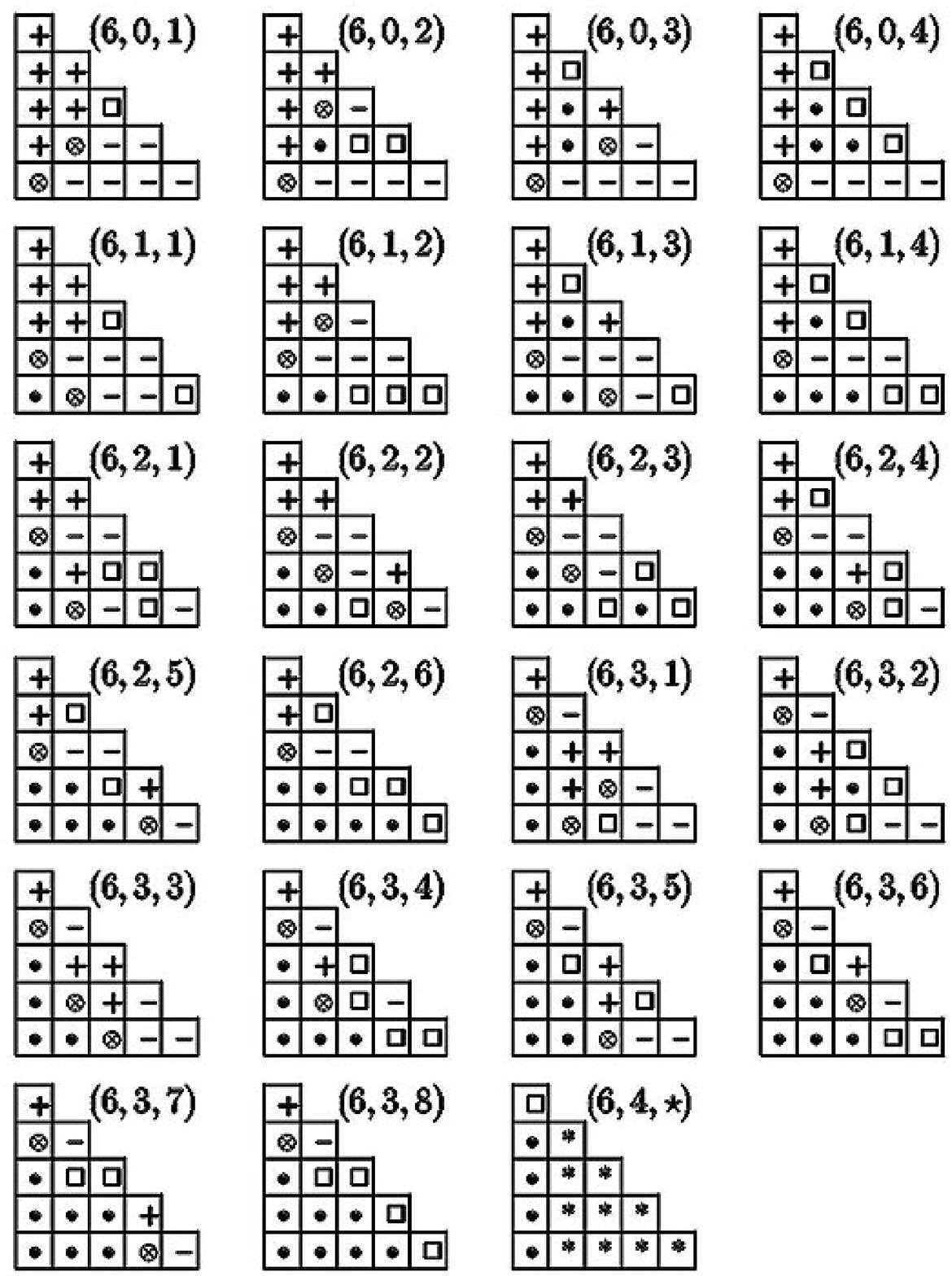}
\newpage
\bf{List of maximal admissible diagrams for $n=7$}
\par\bigskip
\includegraphics[width=16cm]{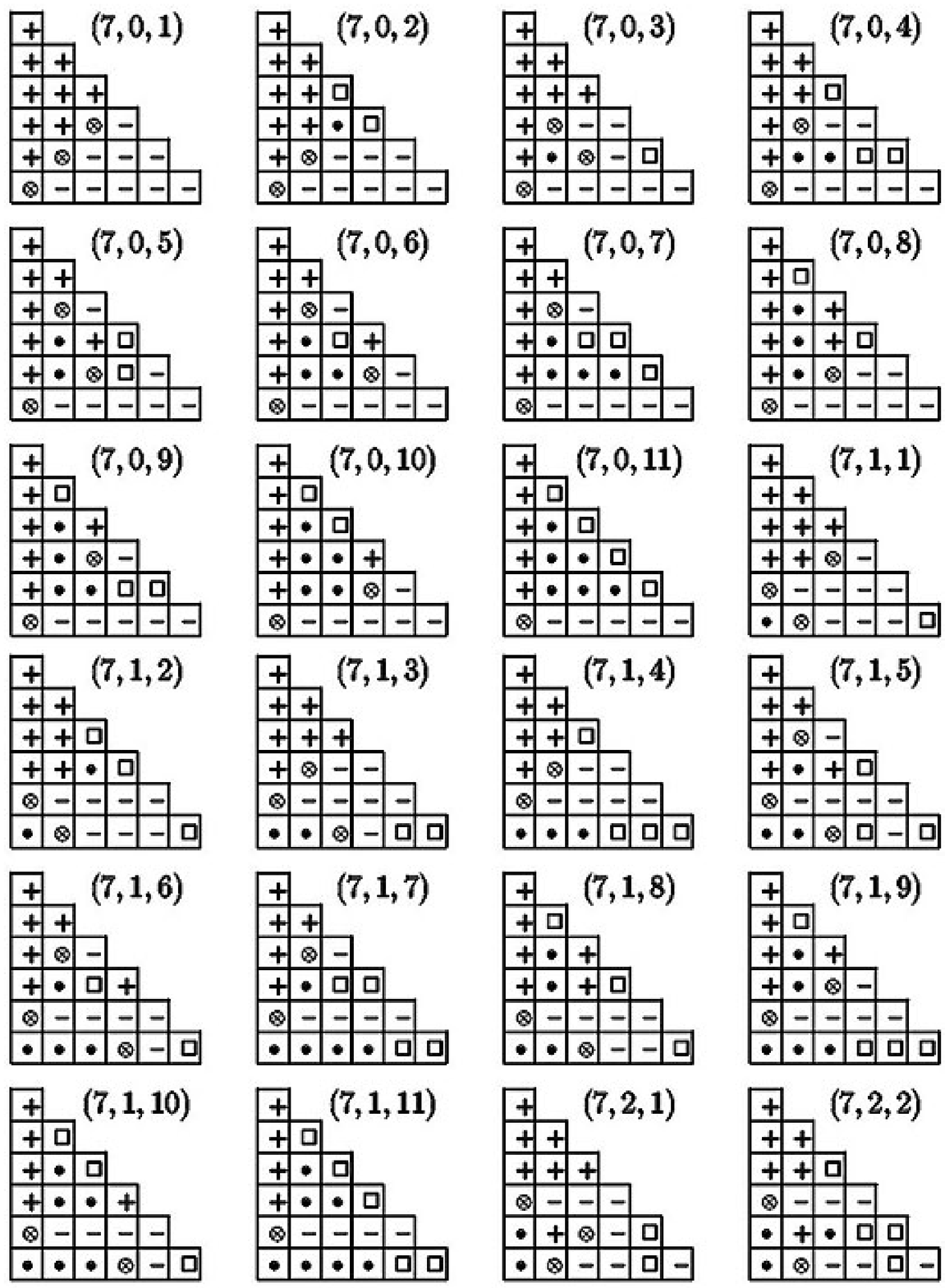}
\newpage
\includegraphics[width=16cm]{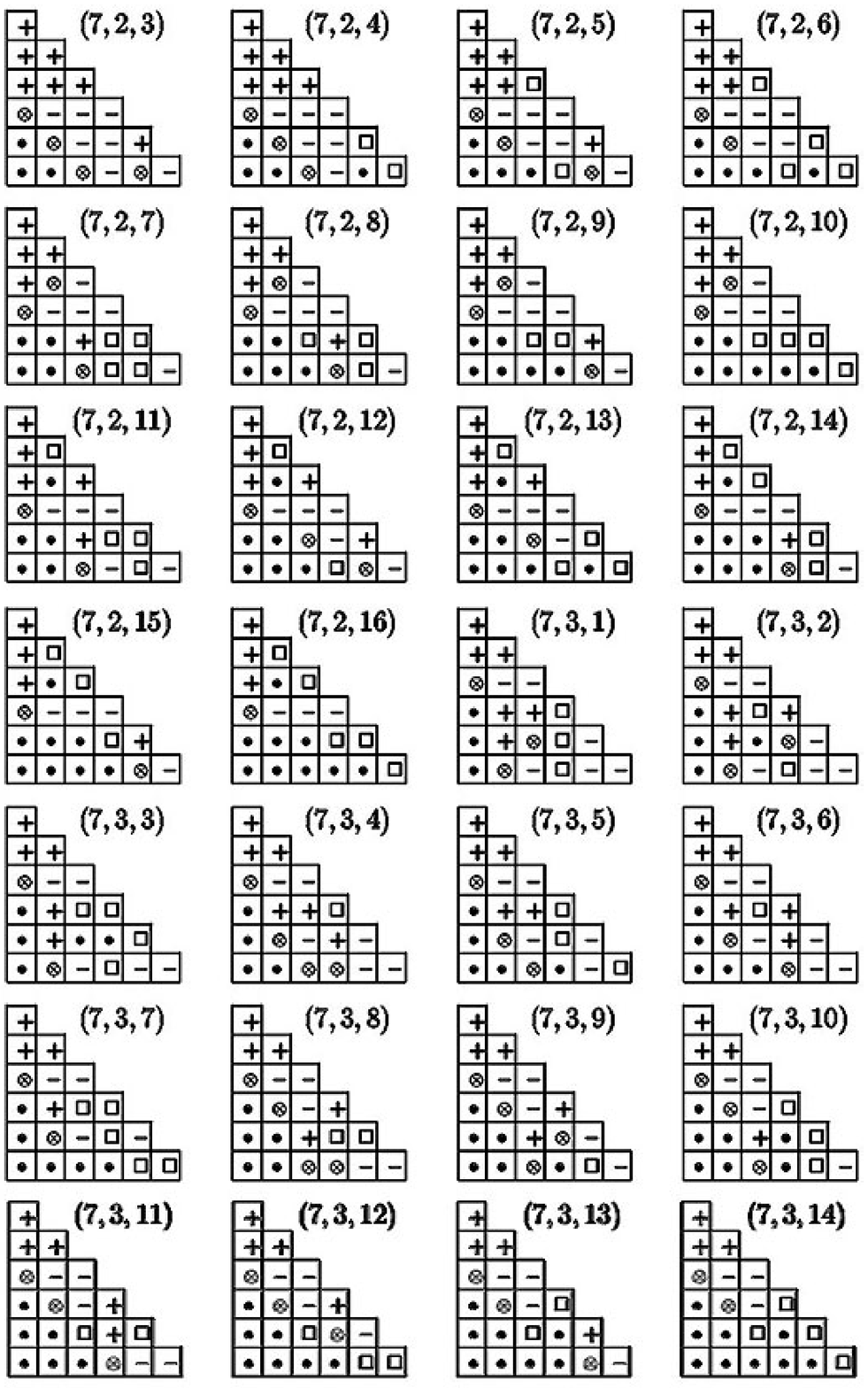}
\newpage
\includegraphics[width=16cm]{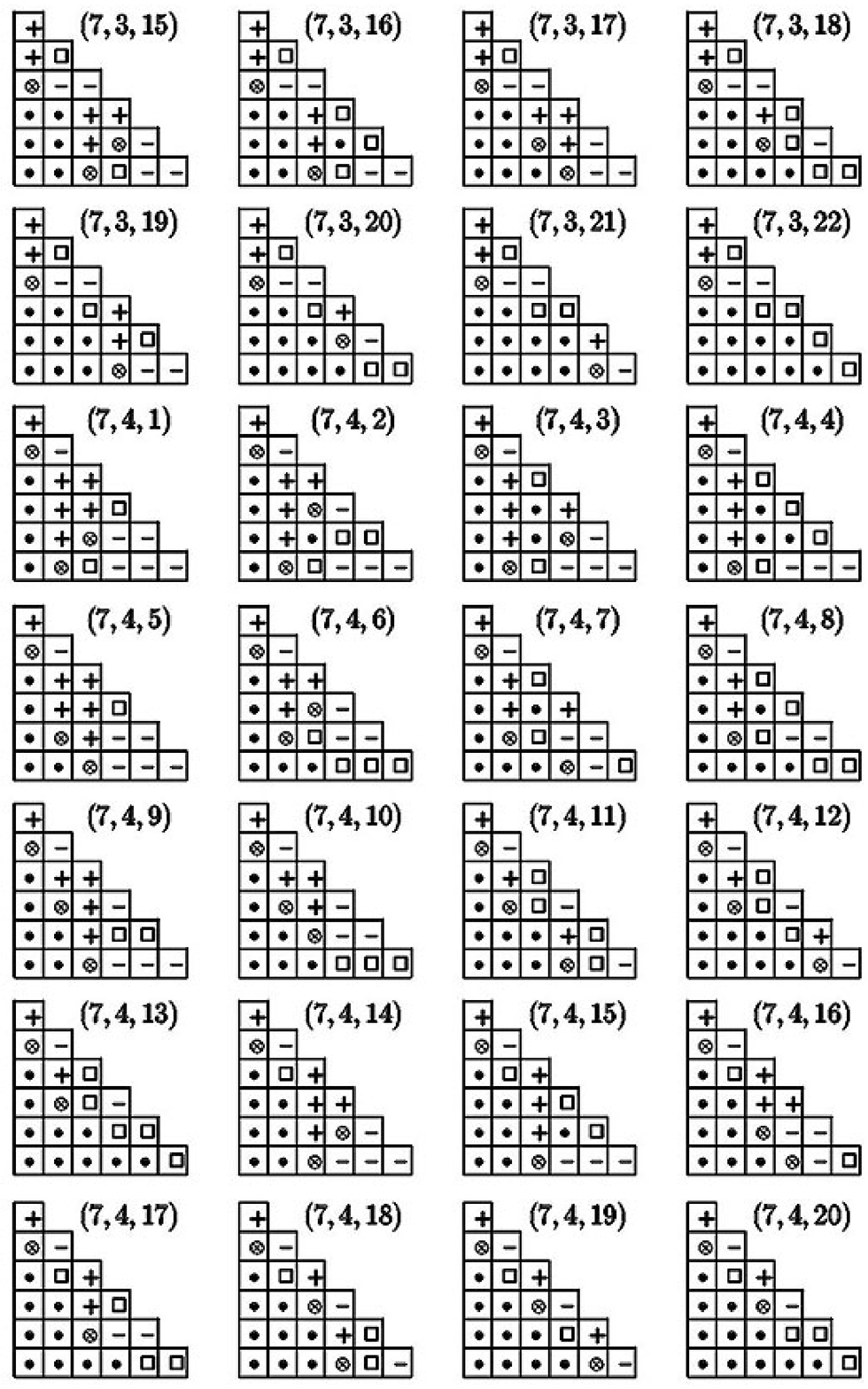}
\newpage
\includegraphics[width=16cm]{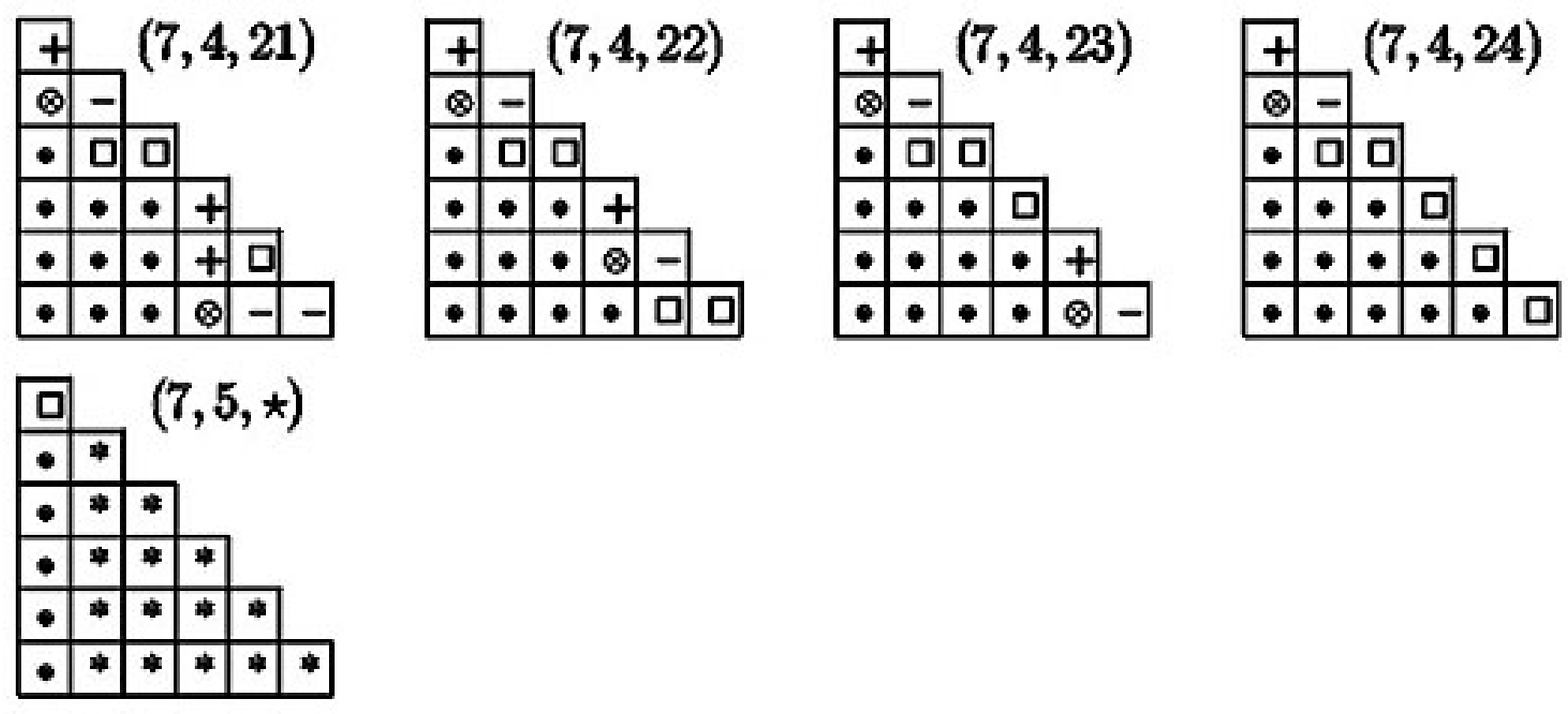}
\par\bigskip\end{center}

\bigskip
Mikhail V. Ignatev\\
Samara State University, Department of algebra and
geometry,\\443011,
ul. Akad. Pavlova, 1, Samara, Russia\\
mihail\_ignatev@mail.ru\par\bigskip

Alexander N. Panov\\
Samara State University, Department of algebra and
geometry,\\443011,
ul. Akad. Pavlova, 1, Samara, Russia\\
apanov@list.ru


\begin{thebibliography}{100}

\bibitem{K-Orb}
Kirillov A.A.. Lectures on the orbit method. Novosibirsk: Nauchnaya
kniga (IDMI), 2002.
\bibitem{K-62} Kirillov A.A. Unitary representations of
nilpotent Lie groups
// Uspekhi matem. nauk. 1962. V.{\bf 17}. P. 57-110.
\bibitem{K-Var}
Kirillov A.A. Variations on the Triangular Theme // Amer. Math. Soc.
Trans.(2). 1995. V.{\bf 169}. P. 43-72.
\bibitem{K-M}
Kirillov A.A., Melnikov A. On a Remarkable Sequences of Polinomials
// in Collection SMF Seminaires et Congres, Algebre non
commutative, groupes quantiques et invariants, J.Alev,
G.Cauchon(editeurs). 1995. \textnumero 2. P. 35-42.
\bibitem{C1}
Andre Carlos A.M. Basic Characters of the Unitriangular Group //
Journal of Algebra. 1995. V.{ \bf 175}. P.287-319.
\bibitem{C2}
Andre Carlos A.M. Basic Sums of Coadjoint Orbits of the
Unitriangular Group // Journal of Algebra. 1995. V.{ \bf 176}.
P.959-1000.
\bibitem{C3}
Andre Carlos A.M. The Basic Character Table of the Unitriangular
Group // Journal of Algebra. 2001. V.{ \bf 241}. P.437-471.
\bibitem{G-Sh}
Gekhtman M.I. Shapiro M.Z. Noncommutative and commutative
Integrability of Generic Toda Flows in Simple Lie Algebras // Comm.
on Pure and Applied Math. 1999. V.{\bf L11}, P.0053-0084.
\bibitem{Bu}
Bourbaki N. Lie groups and Lie algebras. Chapters 4,5,6. Moscow:
Mir, 1972
\bibitem{Dix}
Dixmier J. Universal enveloping algebras. Moscow: Mir, 1978.
\bibitem{Kraft}
Kraft H. Geometric methods in invariant theory. IO NFMI, 2000.

\end{thebibliography}
\end{document}